\documentclass{elsarticle}

\begin{document}

\begin{frontmatter}

\title{Semi-implicit second order schemes for numerical solution of level set advection equation on Cartesian grids}

\author[mymainaddress]{Peter Frolkovi\v{c}\fnref{myfootnote}}
\author[mymainaddress,mysecondaryaddress]{Karol Mikula\fnref{myfootnote}}
\address[mymainaddress]{Slovak University of Technology in Bratislava, Faculty of Civil Engineering, Department of Mathematics and Descriptive Geometry, Radlinsk\'eho 11, 810 05 Bratislava, Slovak Republic}
\ead{peter.frolkovic@stuba.sk, karol.mikula@stuba.sk}
\fntext[myfootnote]{This work was supported by grants VEGA 1/0728/15, VEGA 1/0608/15 and APVV-15-0522.}

\address[mysecondaryaddress]{Algoritmy:SK, \v{S}ulekov\'a 6, 81106 Bratislava, Slovak Republic}

\begin{abstract}
A new parametric class of semi-implicit numerical schemes for a level set advection equation on Cartesian grids is derived and analyzed. An accuracy and a stability study is provided for a linear advection equation with a variable velocity using partial Lax-Wendroff procedure and numerical von Neumann stability analysis. The obtained semi-implicit $\kappa$-scheme is $2^{nd}$ order accurate in space and time in any dimensional case when using a dimension by dimension extension of the one-dimensional scheme that is not the case for analogous fully explicit or fully implicit $\kappa$-schemes. A further improvement is obtained by using so-called Corner Transport Upwind extension in two-dimensional case. The extended semi-implicit $\kappa$-scheme with a specific (velocity dependent) value of $\kappa$ is $3^{rd}$ order accurate in space and time for a constant advection velocity, and it is unconditional stable according to the numerical von Neumann stability analysis for the linear advection equation in general.
\end{abstract}


\end{frontmatter}


\section{Introduction}

In this work we derive a new class of semi-implicit $2^{nd}$ order schemes for numerical solutions of a representative linear advection equation
\begin{equation}
\nonumber
\partial_t u(x,t) + \vec{V} \cdot \nabla u(x,t) = 0 \,, \quad u(x,0) = u^0(x)
\end{equation}
with a variable velocity $\vec{V}=\vec{V}(x)$.
We are interested in level set methods \cite{set99,osh02} when this equation is used to track implicitly given interfaces, and when discontinuous profiles in the solution are not expected in general. The implicit tracking of interfaces can be found in any front propagation problems solved by level set methods, see, e.g., \cite{set99,osh02} and the references there. A typical application is a two-phase flow of immiscible fluids where an interface between the phases must be tracked to distinguish the different physical properties of fluids \cite{sus94,fed99,low06,Herrmann2008a,gro11,Kees2011,Wang2012,flw16}. Furthermore we mention a tracking of fire front in forests \cite{bal12,fmu14}, and a tracking of water table for groundwater flows \cite{her06,fro12}.

We consider Cartesian grids that are often applied in the context of level set methods \cite{sus94,set99,osh02,fw09,fmu14}. 
We consider the linear advection equation on Cartesian grids also as a starting point for a study of more complex equations like a nonlinear advection equation for a motion in normal direction \cite{set99,osh03,fm07,mo10,fmu14} and computations on unstructured grids \cite{fm07,flw16,hmf17}. We are interested here in Eulerian type of numerical schemes of a finite difference form when a stencil of the scheme does not move in time like in Lagrangian type of numerical schemes \cite{fro02,chen2017}. Furthermore we restrict ourselves to the schemes using an implicit or a semi-implicit time discretization with a purpose of favorable stability properties when compared to the schemes of Eulerian type using a fully explicit time discretization.

The fully explicit schemes are standard numerical tool in level set methods for the solution of the linear advection equation \cite{sus94,set99,osh02,xu06,for05, fm07, Min2007,Herrmann2008a,Phongthanapanich2010,gro11,Kees2011,Kim2011}. The main advantage is their simplicity as the numerical solution, once the scheme is constructed, can be obtained directly without solving any algebraic system. 
On the other hand the well-known restriction of fully explicit schemes with fixed stencils is a CFL stability condition on the choice of time steps that depends, among other, on a length of grid steps. 

Although the CFL restriction is not considered as a disadvantage in general, it can be critical for applications with irregular computational domains for which the boundaries are treated implicitly like in Cartesian cut cell methods \cite{lev88b}, immersed interface methods \cite{lev94,li06,xu06,fmu14}, ghost fluid methods \cite{fed99,fed00} and similar. 
In the quoted methods the presence of arbitrary small cut cells can give locally an arbitrary small grid size that results in an unrealistic CFL restriction if no modifications of the numerical scheme is provided.

Recently some publications \cite{mo10,mou14,fmu14} have been dealing with semi-implicit finite volume schemes for a general advection equation. The main idea is that the 
implicit time discretization is used only for the values of numerical solution at inflow boundaries of computational cells.
The semi-implicit schemes can be advantageous when solving the advection equation on implicitly given computational domain as it appears, e.g., when constructing a so-called ``extension'' velocity in level set methods, see \cite{zao96,ada99}. This approach is used in \cite{fmu14} where the linear advection equation is solved by a particular semi-implicit method on a time dependent domain given by positions of a fire front in a forest, and where no cut-cell problem occurs in numerical simulations.

Although some analysis is provided in \cite{mo10,mou14,fmu14} the particular semi-implicit schemes are derived ad hoc. 
In this work we present a unified representation of such semi-implicit schemes using a novel approach of partial Lax-Wendroff procedure and study their accuracy and stability properties.
The Lax-Wendroff \cite{lax60} (or Cauchy-Kowalevski \cite{tor09}) procedure in its full form replaces the time derivatives of the solution in Taylor series by the space derivatives \cite{lev02,tor09}. This procedure is used in a derivation of high order ADER (Arbitrary DERivatives) schemes that are applied to a variety of applications, see, e.g., \cite{tor09} and the references there. In our approach we apply the steps of Lax-Wendroff procedure only partially by allowing the mixed time-space derivatives of the solution in Taylor series.

We use this procedure with an approach of fully explicit $\kappa$-scheme \cite{lee77,lee85,wes01} that includes as particular cases some popular numerical schemes like Lax-Wendroff and Fromm scheme \cite{wes01,lev02} or QUICKEST scheme \cite{leo79,wes01}. The general formulation of the semi-implicit $\kappa$-scheme gives us an opportunity to use special choices of the parameter $\kappa$ to improve the accuracy and the stability of the scheme in special cases, and to adapt the scheme near boundaries. 

To show some advantages of the partial Lax-Wendroff procedure with respect to the full procedure, we compare the semi-implicit $\kappa$-scheme with an analogous fully implicit $\kappa$-scheme derived in this paper using the full Lax-Wendroff procedure. We study the stability conditions of all presented schemes using von Neumann stability analysis \cite{tre96,wes01,hir07} realized in a numerical way as suggested in \cite{bil97a,bil97b}. 
 
The semi-implicit $\kappa$-scheme is unconditionally stable in the one-dimensional case for all relevant values of $\kappa$ that is not the case for the fully implicit $\kappa$-scheme. We show that this property can be used for the immersed interface methods when boundary conditions are defined on an implicitly given boundary of computational domain.
Furthermore we derive a novel particular variant of the semi-implicit $\kappa$-scheme by defining a variable (velocity dependent) value of the parameter $\kappa$. The scheme is $3^{rd}$ order accurate in space and time for a constant velocity in 1D. 

Opposite to the fully implicit $\kappa$-scheme (and also the fully explicit $\kappa$-scheme), the semi-implicit $\kappa$-scheme remains $2^{nd}$ order accurate in space and time in several dimensions when using a standard dimension by dimension extension of 1D scheme on Cartesian grids. Unfortunately, this extension of the semi-implicit $\kappa$-scheme in several dimensions is conditionally stable in general.

To improve the stability of two-dimensional semi-implicit $\kappa$-scheme we apply the idea of Corner Transport Upwind (CTU) extension \cite{col90,lev02} by adding an additional discretization term to the scheme. The main result is a novel scheme with the velocity dependent value of $\kappa$ using the CTU extension that is unconditionally stable according to the numerical von Neumann stability analysis. Moreover the scheme is $3^{rd}$ order accurate in the case of constant velocity. For several representative numerical experiments this variant of the semi-implicit $\kappa$-scheme gives the most accurate results among other considered choices of $\kappa$.

The paper is organized as follows.  In section \ref{sec-1d} we begin with the one-dimensional case where the fully implicit and the semi-implicit $\kappa$-schemes are derived. In 
Section \ref{sec-2d} we discuss the properties of semi-implicit $\kappa$-scheme in several dimensions when obtained by the dimension by dimension extension. Furthermore the Corner Transport Upwind extension of the scheme and the treatment of boundary conditions on implicitly given boundary are described. In section \ref{sec-num} several numerical experiments are presented that involve examples on an implicitly given computational domain, an example with largely varying velocity, and two standard benchmark examples for tracking of interfaces. Finally we conclude the results in section \ref{sec-con}.


\section{One dimensional case}
\label{sec-1d}

We begin with a brief derivation of numerical schemes for a one-dimensional advection equation written as
\begin{equation}
\partial_t u(x,t) + V(x) \partial_x u(x,t) = 0 \,, \quad u(x,0) = u^0(x) \,, \quad x \in (0,L) \,, \quad t \ge 0 \, . 
\label{adveq1d}
\end{equation}
Let $x_i$ be points of a uniform grid such that $x_i - x_{i-1} \equiv h$ with $h=L/M$ and $i=0,1,\ldots,M$.
Furthermore let $\tau>0$ be a given time step and $t^n=n \tau$, $n=0,1,\ldots$ . We use a standard indexing of discrete values like $V_i=V(x_i)$ and so on.

Our aim is to find approximate values $U_i^n$ such that $U_i^n \approx u_i^n$ where $u_i^n := u(x_i,t^n)$ 
and $U_i^0=u^0_i=u^0(x_i)$. 
If $V_0 > 0$ and/or $V_M<0$ then the values $u_0^n$ and/or $u_M^n$ shall be prescribed by Dirichlet boundary conditions. If $V_0 = 0$ and/or $V_M=0$ then $U_0^n=u_0^0$ and/or $U_M^n=u_M^0$. If $V_0 < 0$ and/or $V_M>0$ then we require no boundary conditions, and we use a linear extrapolation by defining  auxiliary values $U_{-1}^n=2 U_0^n-U_1^n$ and $U_{M+1}^n=2 U_M^n-U_{M-1}^n$.

We are interested in finite difference methods using a stencil with at most two neighboring values in each direction, namely
\begin{equation}
\label{all1d}
U_i^{n+1} + \sum_{k=-2}^{2} \alpha_{i k} U_{i+k}^{n+1} = U_i^n + \sum_{k=-2}^{2} \beta_{i k} U_{i+k}^{n} \,.
\end{equation}
A fully explicit form of (\ref{all1d}) is given by $\alpha_{i k} \equiv 0$, analogously $\beta_{ i k} \equiv 0$ in the case of fully implicit form.
Our aim is to derive semi-implicit schemes that have in general three consecutive nonzero values of coefficients $\alpha_{i k}$ and $\beta_{i k}$ in (\ref{all1d}). 

To check an order of accuracy for any particular scheme of the form (\ref{all1d}) we consider its truncation error that is obtained by replacing all numerical values $U_{i+k}^{n+1}$ and $U_{i+k}^{n}$  in (\ref{all1d}) with the exact values $u_{i+k}^{n+1}$ and $u_{i+k}^{n}$ that themselves are then expressed with Taylor series, see some standard textbooks on numerical analysis, e.g., \cite{hir07}.

To derive a stability condition of particular numerical scheme of he form (\ref{all1d}), we use an approach of von Neumann stability analysis, see, e.g., \cite{tre96,wes01,hir07}. To do so one introduces a grid function $\epsilon_i^n=\epsilon(x_i,t^n)$ defined by  
\begin{equation}
\label{eps1d}
\epsilon(x,t)=\exp(-\lambda t) \exp(\imath x) \,, \quad x \in R \,, \,\, t \ge 0 \,,
\end{equation}
where $\imath$ is the imaginary number, and the parameter $\lambda$ shall be found. 
The values $\epsilon_i^n$ are supposed to fulfill the numerical scheme (\ref{all1d}). Using  relations
\begin{equation}
\label{trivial}
\epsilon_{i\pm k}^n =\exp(\pm \imath k h) \epsilon_i^n \,, \quad
\epsilon_i^{n+1}=S \epsilon_i^n \,, \,\, S:=\exp(-\lambda \tau) \,, 
\end{equation}
where $S$ denotes an amplification factor, the {von Neumann stability analysis} is realized by searching for conditions under which one has $|S| \le 1$ for all $h \in (-\pi,\pi)$ . Using (\ref{trivial}) in (\ref{all1d}) one obtains
\begin{equation}
\label{allstab}
S = \left(1 + \sum_{k=-2}^{2} \beta_{i k} \exp(\imath k h) \right) \left(1+ \sum_{k=-2}^{2} \alpha_{i k}  \exp(\imath k h) \right)^{-1} \,.
\end{equation}

Although the stability conditions for $|S| \le 1$ from (\ref{allstab}) can be found using analytical methods  for some schemes \cite{tre96,wes01,hir07}, we apply an approach proposed and used in \cite{bil97a,bil97b} where such condition is found numerically. One approach is to compute the values $|S|$  for very large number of discrete values of $h$ and the input parameters of particular numerical scheme \cite{bil97a,bil97b}. We apply numerical optimization algorithms available in Mathematica$^{\tiny{\textregistered}}$ \cite{mat16} to search for local maxima of $|S|$ using a very large number of different initial guesses. Such numerical von Neumann stability analysis is applied to all numerical schemes studied in this paper including nontrivial two-dimensional schemes later.

A so-called fully explicit $\kappa$-scheme \cite{lee77,lee85,wes01} of the form (\ref{all1d}) for a variable advection velocity can be found in \cite{fm17}. 
We describe now the derivation of $2^{nd}$ order accurate fully implicit $\kappa$-scheme to solve (\ref{adveq1d}). The accuracy is obtained in space and time that we do not  emphasize furthermore in our description.

In what follows we use shorter notations for the exact values of derivatives by $\partial_t u_i^* := \partial_t u(x_i,t^*)$ for $*=n$ or $*=n+1$ and so on. An analogous notation with the capital letter $U$ is reserved for numerical approximations of derivatives. An important role in our derivation will play the following parametric class of approximations $\partial_x^{\kappa} U_i^* \approx \partial_x u_i^*$ 
\begin{eqnarray}
\label{kappagrad}
2 \, \partial_x^{\kappa} U_{i}^* := (1-\kappa) \, \partial_x^- U_i^* + (1+\kappa) \, \partial_x^+ U_i^* \,, 
\end{eqnarray}
where
$$
h \, \partial_x^- U_i^* := U_i^*-U_{i-1}^* \,, \quad h \, \partial_x^+ U_i^* := U_{i+1}^*-U_i^* \,,
$$
and the parameter $\kappa$ in (\ref{kappagrad})  is free to choose. A natural choice $\kappa \in [-1,1]$ gives a convex combination of the standard one-sided finite difference approximations.

Another important tool to derive the scheme is Lax-Wendroff procedure \cite{lax60}, also quoted as Cauchy-Kowa\-lewski procedure \cite{tor09}, that consists of replacing all time derivatives of $u$ by the space derivatives of $u$ using the equation (\ref{adveq1d}). We write it in the form
\begin{eqnarray}
\label{ck01}
\partial_t u_i^{n+1} = - V_i \partial_x u_i^{n+1} \,,\\
\label{ck02}
\partial_{tx} u_i^{n+1} = - \partial_x (V \partial_{x} u)_i^{n+1} \,,\\
\label{ck11}
\partial_{tt} u_i^{n+1} = - V_i \partial_{tx} u_i^{n+1} = V_i \partial_x (V \partial_{x} u)_i^{n+1} \,.
\end{eqnarray}

Now using the Taylor series in a backward manner
\begin{equation}
\label{imt}
u_i^{n} = u_i^{n+1} + \sum_{m=1}^p \frac{(-1)^m}{m!} \tau^m \partial^m_t u_i^{n+1}  + \mathcal{O}(\tau^{p+1})
\end{equation}
and the Lax-Wendroff procedure (\ref{ck01}) - (\ref{ck11}) with (\ref{imt}) for $p=2$ one obtains
\begin{eqnarray}
\label{accur2i}
u_i^{n} = u_i^{n+1} + \tau V_i \partial_x u_i^{n+1} + 0.5 \tau^2 V_i \partial_{x} (V \partial_x u)_i^{n+1}
+ \mathcal{O}(\tau^3) \,.
\end{eqnarray}

The fully implicit $\kappa$-scheme to solve (\ref{adveq1d}) is now obtained by applying proper (upwinded) finite difference approximations in (\ref{accur2i}). To reach a truncation error corresponding to $2^{nd}$ order accurate schemes, the term multiplied by $\tau$ in (\ref{accur2i}) must be approximated by a $2^{nd}$ order accurate approximation in space, while for the term multiplied by $0.5 \tau^2$ a first order accurate approximation in space is sufficient.
Having this in mind we apply in (\ref{accur2i}) the following upwinded approximations
\begin{eqnarray}
\label{upwgrad1d}
V_i \partial_x u_i^{n+1} \approx
[V_i]^+\partial_x^- (U_i^{n+1} + 0.5 h \partial^{\kappa}_x U_i^{n+1}) +[V_i]^- \partial_x^+ (U_i^{n+1} - 0.5 h \partial^{\kappa}_x U_i^{n+1})  \,, \\[1ex]\label{dxx}
V_i \partial_{x} (V \partial_x u)_i^{n+1} \approx \left( [V_i]^+ \partial^-_x + [V_i]^- \partial^+_x\right) (V_i \partial_x^{\kappa} U_i^{n+1})
\end{eqnarray}
where $[V]^+=\max\{0,V\}$ and $[V]^-=\min\{0,V\}$. 
Using (\ref{upwgrad1d}) - (\ref{dxx}) in (\ref{accur2i}), the fully implicit $\kappa$-scheme is obtained
\begin{eqnarray}
\label{impl1d}
U_i^{n+1} +  \tau [V_i]^+ \partial_x^- \left( U_i^{n+1} + 0.5 (h + \tau  V_i) \partial_x^{\kappa} U_i^{n+1})\right) + \\[1ex] 
\nonumber
+ \left. \tau [V_i]^- \partial_x^+ \left( U_i^{n+1} - 0.5 (h - \tau  V_i) \partial_x^{\kappa} U_i^{n+1})\right) \right.  = U_i^{n}\,.
\end{eqnarray}
Note that due to the presence of the term $\partial_x^- \left(V_i \partial_x^{\kappa} U_i^{n+1}\right) = \left(V_i \partial_x^{\kappa} U_i^{n+1} - \right. $ $\left.V_{i-1} \partial_x^{\kappa} U_{i-1}^{n+1}\right) / h$, and analogously for $\partial_x^+ \left(V_i \partial_x^{\kappa} U_i^{n+1}\right)$, the scheme (\ref{impl1d}) uses two discrete values of velocity.

Similarly to the fully explicit $\kappa$-scheme \cite{fw09,fm17}, the scheme (\ref{impl1d}) is $2^{nd}$ order accurate for an arbitrary value of $\kappa$.
Analogously to the fully explicit QUICKEST scheme \cite{leo79,wes01}, one can prove that the choice
\begin{equation}
\label{im3rd}
\kappa = \hbox{sign}(\mathcal{C}) ( 1 + 2 |\mathcal{C}| ) / 3 \,, \quad {\mathcal C} := \tau V/h
\end{equation}
gives the $3^{rd}$ order accurate scheme in space and time in the case of constant velocity $V$.

We discuss now the von Neumann stability analysis \cite{tre96,wes01,hir07} for (\ref{impl1d}) that is realized for locally frozen values of the velocity $V(x) \equiv V_i$. We denote the (signed) grid Courant numbers by
$$
{\mathcal C}_i := \tau V_i / h \,.
$$
The numerical von Neumann stability analysis shows that for $\kappa \le 0$ the fully implicit $\kappa$-scheme is unconditionally stable for $\mathcal{C}_i\ge 0$. In the case $\mathcal{C}_i \le 0$ the unconditional stability is obtained for $\kappa \ge 0$.
These stability conditions are more favorable when compared to conditions of the fully explicit $\kappa$-scheme. The price to pay is that a system of linear algebraic equations (\ref{impl1d}) has to be solved in each time step to obtain the values $U_i^{n+1}$. 

Unfortunately, other interesting choices of $\kappa$ give only restrictive stability conditions. For instance 
the value $\kappa=1/3$ gives the stability condition $0 \le \mathcal{C}_i \le 2$, the $3^{rd}$ order accurate choice (\ref{im3rd}) gives the condition $|\mathcal{C}_i| \le 0.5$. Moreover, the choice $\kappa=\hbox{sign}(V_i)$ results in an unstable scheme for $\mathcal{C}_i \in [-1,1]$. 

Next we present the semi-implicit variant of $\kappa$-scheme.
The main idea is to apply the partial Lax-Wendroff procedure by skipping the replacement of $\partial_{tx} u$ in (\ref{ck11}). Using (\ref{ck01}) and (\ref{ck02}) with (\ref{imt}) for $p=2$ we obtain
\begin{eqnarray}
\label{accur4}
u_i^{n} =  u_i^{n+1} + \tau V_i \partial_x u_i^{n+1} - 0.5 \tau^2 V_i \partial_{t x}  u_i^{n+1}
+ \mathcal{O}(\tau^3) \,.
\end{eqnarray}
Now using the approximation (\ref{upwgrad1d}) in (\ref{accur4})  for the term containing $\partial_x u_i^{n+1}$ and the following approximation for the second term
$$
\tau \, \partial_{tx} u_i^{n+1} \approx \tau \, \partial_t^- \partial_x^{\kappa} U_i^{n+1} = \partial_x^{\kappa} U_i^{n+1} - \partial_x^{\kappa} U_i^{n} \,,
$$
one obtains 
\begin{eqnarray}
\nonumber
U_i^{n} = U_i^{n+1} +  \tau [V_i]^+ \left(\partial_x^- U_i^{n+1} + 0.5 h \partial_x^-  \partial_x^{\kappa} U_{i}^{n+1} \right) + \\[1ex]
\nonumber
+ \left. \tau [V_i]^- \left(\partial_x^+ U_i^{n+1} - 0.5 h \partial_x^+ \partial_x^{\kappa} U_{i}^{n+1}  \right) \right. 
-0.5 \tau V_i \left( \partial_x^{\kappa} U_i^{n+1} - \partial_x^{\kappa} U_i^{n} \right)\,.
\end{eqnarray}
After simple algebraic manipulations the semi-implicit $\kappa$-scheme can be written in the form
\begin{eqnarray}
\label{si1d}
U_i^{n+1} +  \tau V_i \left(\partial_x^{\mp} U_i^{n+1} - 0.5 \partial_x^{\kappa} U_{i\mp 1}^{n+1} \right) 
= \left.
U_i^{n} -   0.5 \tau V_i \partial_x^{\kappa} U_{i}^{n} \right. \,,
\end{eqnarray}
where one has to replace $\mp$ in $\partial_x^{\mp}$ and in $i\mp 1$ with the opposite signs with respect to $\hbox{sign}(V_i)$, i.e. with $-$ if $V_{i}>0$ and with $+$ if $V_{i}<0$. 
Note that opposite to the fully implicit $\kappa$-scheme (\ref{impl1d}), the semi-implicit variant (\ref{si1d}) requires locally only the single value $V_i$.

Looking at the truncation error of  (\ref{si1d}) the scheme is $2^{nd}$ order accurate for an arbitrary value of $\kappa$. The choice 
\begin{equation}
\label{kappasi}
\kappa = \hbox{sign}(\mathcal{C}_i) (1- |\mathcal{C}_i|)/3
\end{equation}
gives in the case of a constant velocity, i.e. ${\mathcal C}_i \equiv {\mathcal C}$, the $3^{rd}$ order accurate scheme in space and time.  As we show later in numerical experiments the choice (\ref{kappasi}) is advantageous also for the variable velocity case.

The most important property of semi-implicit $\kappa$-scheme is its stability condition. Opposite to the fully implicit $\kappa$-scheme (\ref{impl1d}), the semi-implicit $\kappa$-scheme (\ref{si1d}) exploits only the single value $V_i$, so the stability analysis is valid without locally freezing the velocity. The numerical von Neumann stability analysis implies unconditional stability for arbitrary $\kappa \le 1$ if $V_i>0$ and for $\kappa \ge -1$ if $V_i<0$ that is a clear improvement with respect to the fully implicit $\kappa$-scheme.

We present three particular variants of (\ref{si1d}) to present them in a more clear way. In fact, the following choices are our suggestions. As we describe later, the schemes can be used in several dimensions by the standard dimension by dimension extension for Cartesian grids. 

Firstly, for the choice $\kappa = \hbox{sign}(V_i)$ the scheme (\ref{si1d}) takes the form
\begin{eqnarray}
\label{siLW}
U_i^{n+1} +  0.5 \, |\mathcal{C}_i| ( U_i^{n+1} - U_{i\mp 1}^{n+1} ) = U_i^{n} -  0.5 \, |\mathcal{C}_i| ( U_{i\pm 1}^{n} - U_{i}^n ) \,,
\end{eqnarray}
where the sign in $\mp$ is chosen opposite to the sign of $\mathcal{C}_i$ and the sign in $\pm$ is identical to the sign of $\mathcal{C}_i$. The scheme (\ref{siLW}) is introduced in \cite{mo11} in a finite volume context in several dimensions, see also \cite{mo10,mou14}. The scheme has the smallest stencil in the implicit part among all particular variants of $\kappa$-scheme, therefore we recommend to use it for the grid nodes next to inflow boundaries.  The amplification factor $S$ of this scheme equal to $1$ everywhere, so any oscillations in numerical solutions are not damped, and the method may require a limiting (or a stabilization \cite{mou14}) procedure in general, see also related numerical experiments. 

Secondly, for the choice $\kappa \equiv 0$ the scheme (\ref{si1d}) turns to
\begin{eqnarray}
\label{siF}
U_i^{n+1} +  0.25 \, |\mathcal{C}_i| ( 3 U_i^{n+1} - 4 U_{i\mp 1}^{n+1} +U_{i\mp 2}^{n+1} ) = U_i^{n} -  0.25 \, \mathcal{C}_i ( U_{i+1}^{n} - U_{i-1}^n ) \,.
\end{eqnarray}
The scheme (\ref{siF}) is introduced in \cite{fmu14} in the finite volume context. The scheme is based on the central difference for the approximation of $\partial_x u$ in (\ref{kappagrad}) , and it gives for the examples presented in the section on numerical experiments the most accurate results among all considered constant values of $\kappa$ parameters. As we discuss later the scheme is not unconditionally stable in several dimensions when using the dimension by dimension extension, but it has a much less restrictive stability condition than the analogous fully explicit scheme.

Finally, for the velocity dependent choice (\ref{kappasi}) the scheme (\ref{si1d}) gets the form
\begin{eqnarray}
\label{siQ}
U_i^{n+1} +  \frac{|\mathcal{C}_i|}{6} \left( 4 U_i^{n+1} - 5 U_{i\mp 1}^{n+1} + U_{i\mp 2}^{n+1} \right) + \frac{\mathcal{C}_i^2}{12} \left( U_i^{n+1} - 2 U_{i\mp 1}^{n+1} + U_{i\mp 2}^{n+1} \right) = \\[1ex] 
\nonumber
U_i^{n} -  \frac{|\mathcal{C}_i|}{6} \left( 2 U_{i\pm 1}^{n} - U_i^n - U_{i\mp 1}^n \right) +  \frac{\mathcal{C}_i^2}{12} \left( U_{i+1}^{n} -2 U_i^n + U_{i-1}^n \right) ,
\end{eqnarray}
where the signs $\pm$ and $\mp$ are chosen as in (\ref{siLW}). The scheme gives the most accurate results for the chosen examples in the section on numerical experiments among all considered variants of $\kappa$-scheme. As we describe in the next section using the Corner Transport Upwind extension the scheme is unconditionally stable in two-dimensional case when using its dimension by dimension extension. 

\section{Two-dimensional case}
\label{sec-2d}

The representative two-dimensional advection equation takes the form
\begin{equation}
\partial_t u(x,y,t) + \vec{V}(x,y) \cdot \nabla u(x,y,t) = 0 \,, \quad
u(x,y,0)=u^0(x,y) \,,
\label{adveq}
\end{equation}
where $\vec{V}=(V(x,y),W(x,y))$.
We restrict ourselves to Cartesian grids that are obtained from the uniform one-dimensional grids using the standard dimension by dimension extension. We denote the uniform space discretization step by $h$.  Our aim is to determine the approximate values $U_{ij}^n \approx u_{i j}^{n}$ where  $u_{ij}^n := u(x_i,y_j,t^n)$.

An extension of the general 1D scheme (\ref{all1d}) for a two-dimensional case can be written in the form
\begin{eqnarray}
\label{all2d}
U_{ij}^{n+1} + \\
\nonumber
\sum_{k=-2}^{2} \left( \alpha^x_{i j k} U_{i+k j}^{n+1} + \alpha^y_{i j k} U_{i j+k}^{n+1} \right) = 
U_{ij}^{n} +\sum_{k=-2}^{2} \left( \beta^x_{i j k} U_{i+k j}^{n} + \beta^y_{i j k} U_{i j+k}^{n} \right) ,
\end{eqnarray}
where we adopt a notation of superscripts $x$ and $y$ to relate the coefficients to a particular space variable.

Analogously to (\ref{kappagrad}), the following approximation of gradients $( \partial_x u_{i j}^n,\partial_y u_{i j}^n)$ is used,
\begin{eqnarray}
\label{grad2d}
2 h \, \partial_x^{\kappa} U_{i j}^n =  (1-\kappa^x) \partial_x^- U_{i j}^n + (1+\kappa^x) \partial_x^+ U_{i j}^n \,, \\
\label{grad2dB}
2 h \, \partial_y^{\kappa} U_{i j}^n =  (1-\kappa^y) \partial_y^- U_{i j}^n + (1+\kappa^y) \partial_y^+ U_{i j}^n\,,  
\end{eqnarray}
where $\partial_x^{-}$, $\partial_x^{+}$, $\partial_y^{-}$, and $\partial_y^{+}$ denote the standard finite differences analogously to (\ref{kappagrad}), and the parameters $\kappa^x$ and $\kappa^y$ are free to choose. 

To provide a stability analysis of (\ref{all2d}), one extends the one-dimensional treatment (\ref{eps1d}) - (\ref{allstab}) by using a grid function $\epsilon_{ij}^n=\epsilon(x_i,y_j,t^n)$ 
when the amplification factor $S$ takes the form analogous to (\ref{allstab})
\begin{eqnarray}
\label{allstab2d}
S = \frac{
1 + \sum_{k=-2}^{2} (\beta^x_{i j k} \exp(\imath k x)+\beta^y_{i j k} \exp(\imath k y)) 
}{
1+ \sum_{k=-2}^{2} \alpha^x_{i j k}  \exp(\imath k x) + \sum_{k=-2}^{2} \alpha^y_{i j k}  \exp(\imath k y) 
}
\end{eqnarray}
and $x,y \in (-\pi,\pi)$.

Before discussing some particular numerical schemes of the form (\ref{all2d}) we note that the Lax-Wendroff procedure takes now more involved form than (\ref{ck01}) - (\ref{ck11}) in 1D case, as it contains also mixed derivatives, namely
\begin{eqnarray}
\label{taylor2dA}
\partial_{t} u_{ij}^n = -V_{ij} \partial_{x} u_{ij}^n - W_{ij} \partial_{y} u_{ij}^n \,,  \,\, \partial_{tt} u_{ij}^n = -V_{ij} \partial_{tx} u_{ij}^n - W_{ij} \partial_{ty} u_{ij}^n \,,  \,\, \\[1ex]
\label{taylor2dB}
\partial_{tx} u_{ij}^n = -\partial_x (V \partial_{x} u)_{ij}^n - \partial_x (W \partial_{y} u)_{ij}^n,  
\partial_{ty} u_{ij}^n = -\partial_y (V \partial_{x} u)_{ij}^n - \partial_y (W \partial_{y} u)_{ij}^n.
\end{eqnarray}

Clearly, when using the full Lax-Wendroff procedure to derive a fully implicit variant of (\ref{all2d}), one has to approximate, due to (\ref{taylor2dB}), the mixed spatial derivative of $u$. This can not be done using the stencil prescribed by (\ref{all2d}), consequently the dimension by dimension extension of the fully implicit 1D scheme (\ref{impl1d}) is not $2^{nd}$ order accurate in space and time as already well-known from a literature for the fully explicit 1D schemes, see, e.g., \cite{bil97b,lev02}.

On the other hand, when using the partial Lax-Wendroff procedure (\ref{taylor2dA}) without (\ref{taylor2dB}), the mixed spatial derivative is not  involved. As a consequence the dimension by dimension extension of 1D semi-implicit $\kappa$-scheme (\ref{si1d}) is $2^{nd}$ order accurate for a variable velocity $\vec{V}(x,y)$ and for arbitrary values of $\kappa^x$ and $\kappa^y$ in (\ref{grad2d}) and (\ref{grad2dB}). 

Similarly to (\ref{si1d}) we can write the semi-implicit $\kappa$-scheme in the compact (upwind) way 
\begin{eqnarray}
\nonumber
U_{i j}^{n+1} + 
\tau V_{i j} \left(\partial_x^{\mp} U_{i j}^{n+1} - 0.5 \partial_x^{\kappa} U_{i\mp 1\, j}^{n+1} \right) + 
\tau W_{i j} \left(\partial_y^{\mp} U_{i j}^{n+1} - 0.5 \partial_y^{\kappa} U_{i\, j\mp 1}^{n+1} \right) 
= \\[1ex]
\label{si2d}
U_i^{n} - 0.5 \tau \left(V_{i j} \partial_x^{\kappa} U_{i j}^{n} + W_{i j} \partial_y^{\kappa} U_{i j}^{n} \right) ,
\end{eqnarray}
where one has to replace $\mp$ in $\partial_x^{\mp}$ and in $i\mp 1$ with $-$ if $V_{ij}>0$ and with $+$ if $V_{ij}<0$, compare also with (\ref{si1d}),  and analogously for the cases related to the sign of $W_{ij}$. Note that the scheme (\ref{si2d}) can be easily extended to a higher dimensional case.

It appears that, opposite to 1D case, the semi-implicit $\kappa$-scheme (\ref{si2d}) is conditionally stable in general. It is not easy to characterize the  stability condition for (\ref{si2d}) as the amplification factor $S$ in (\ref{allstab2d}) depends on six free parameters: $x$, $y$, $\kappa^x$, $\kappa^y$ and two (directional and signed) grid Courant numbers
\begin{eqnarray}
\label{cn2d}
\mathcal{C}_{ij}= \frac{\tau V_{ij}}{h} \,, \,\, \mathcal{D}_{i j} = \frac{\tau W_{ij}}{h} \,.
\end{eqnarray}


We give here some details about stability conditions for the variants (\ref{siLW}), (\ref{siF}), and (\ref{siQ}) used in the form (\ref{si2d}).
The choice $\kappa^x=\hbox{sign}(V_{ij}), \kappa^y=\hbox{sign}(W_{ij})$ gives so-called IIOE scheme (Inflow Implicit / Outflow Explicit) published in a finite volume form in \cite{mou14}. The scheme gives $|S|=1$ for all values of $\mathcal{C}_{ij},\mathcal{D}_{ij},x,y$. Consequently it is unconditionally stable, but it does not damp any oscillations in numerical solutions, see related numerical experiments later.

The choice $\kappa^x= \kappa^y=0$ represented by (\ref{siF}) is used in a finite volume form in \cite{fmu14}. The numerical von Neumann stability analysis gives the  stability condition $|S| \le 1$ for $|\mathcal{C}_{ij}| \le 7.396$ and $|\mathcal{D}_{ij}| \le 7.396$ that is significantly less restrictive than in the case of fully explicit schemes. The value $|S|$ can be larger than 1 in general, for instance, the maximal value of $|S|$ is around $1.00013$ and $1.04538$ for the maximal Courant numbers $8$ and $16$, respectively. 

Finally, the variable choice (\ref{kappasi}) of $\kappa$ represented by (\ref{siQ}) is stable for $|\mathcal{C}_{ij}| \le 4$ and $|\mathcal{D}_{ij}| \le 4$. In the next section we extend the semi-implicit $\kappa$-scheme in two-dimensional case to such a form for which the unconditional stability is obtained by numerical von Neumann stability analysis for the dimension by dimension extension of (\ref{siQ}).

\subsection{Corner Transport Upwind extension}

We begin our description with an analysis of the truncation error of (\ref{si2d}) for the case of constant velocity $\vec{V}$. We keep the indexing in the notation $(V_{ij},W_{ij})$ of discrete velocity values.

We consider the Taylor series
\begin{equation}
\label{imt2d}
u_{ij}^{n} = u_{ij}^{n+1} + \sum_{m=1}^p \frac{(-1)^m}{m!} \tau^m \partial^m_t u_{ij}^{n+1}  + \mathcal{O}(\tau^{p+1}) \,,
\end{equation}
for $p=3$ and replace the first and second time derivatives in (\ref{imt2d}) using (\ref{taylor2dA}) at $t^{n+1}$. The third term can be replaced for the constant velocity $(V_{ij},W_{ij})$ by
\begin{eqnarray}
\label{taylor2d3A}
\partial_{ttt} u_{ij}^{n+1} = -V_{ij} \partial_{ttx} u_{ij}^{n+1} - W_{ij} \partial_{tty} u_{ij}^{n+1} \,.
\end{eqnarray}
Furthermore we use  in (\ref{taylor2d3A}) the following relations  that are valid for the constant vector $(V_{ij},W_{ij})$,
\begin{eqnarray}
\label{taylor2d3B}
\partial_{ttx} u_{ij}^{n+1} = -V_{ij} \partial_{txx} u_{ij}^{n+1} - W_{ij}  \partial_{txy} u_{ij}^{n+1} \,,  \,\,\\[1ex]
\label{taylor2d3C}
\partial_{tty} u_{ij}^{n+1} = -V_{ij} \partial_{txy} u_{ij}^{n+1} - W_{ij} \partial_{tyy} u_{ij}^{n+1}\,.
\end{eqnarray}

Now one can show that for the choices analogous to (\ref{kappasi})
\begin{equation}
\kappa^x = \hbox{sign}(\mathcal{C}_{i j}) (1- |\mathcal{C}_{i j}|)/3 \,, \,\,
\kappa^y = \hbox{sign}(\mathcal{D}_{i j}) (1- |\mathcal{D}_{i j}|)/3 \,, 
\label{kappasi2d}
\end{equation}
the spatial derivatives $\partial_{xxx}u_{i j}^{n+1}$ and $\partial_{yyy} u_{i j}^{n+1}$  are canceled in the truncation error analogously to the 1D case. Nevertheless, the following third order derivatives term 
will remain
\begin{equation}
\frac{\tau^2}{12} V_{i j} W_{i j} \left(2 \tau \partial_{txy} u_{i j}^{n+1} - h \partial_{xxy} u_{i j}^{n+1} - h \partial_{xyy} u_{i j}^{n+1} \right) \,.
\label{mixed}
\end{equation}

In what follows we apply an idea of Corner Transport Upwind (CTU) extension \cite{col90} to extend the scheme (\ref{si2d}) in a such way that also the term in (\ref{mixed}) will cancel in the truncation error. We follow \cite{lev02} where it is used to extend the fully explicit schemes of the form (\ref{all2d}).   

The CTU extension adds an additional discretization term to (\ref{si2d}) that contains, additionally to the stencil of (\ref{all2d}), also the corner (diagonal) values of numerical solution. We do it in a such way that the $2^{nd}$ order accuracy of (\ref{si2d}) is preserved, and the scheme (\ref{si2d}) with the CTU extension using the variable choice (\ref{kappasi2d}) of $\kappa$ parameters becomes $3^{rd}$ order accurate in the case of a constant velocity $\vec{V}$.

In fact one can derive a parametric class of the scheme with the CTU extension where its particular variants are different only in their explicit part, and that can be obtained as a convex combination of two representative schemes,
\begin{eqnarray}
\nonumber
U_{i j}^{n+1} + 
\tau V_{i j} \left(\partial_x^{\mp} U_{i j}^{n+1} - 0.5 \partial_x^{\kappa} U_{i\mp 1\, j}^{n+1} \right) + 
\tau W_{i j} \left(\partial_y^{\mp} U_{i j}^{n+1} - 0.5 \partial_y^{\kappa} U_{i\, j\mp 1}^{n+1} \right) + \\[1ex]
\nonumber
|\mathcal{C}_{ij} \mathcal{D}_{ij}|/6 \left( U_{i j}^{n+1}+U_{i\mp1\, j\mp1}^{n+1}-U_{i\mp1\, j}^{n+1}-U_{i\, j\mp1}^{n+1} \right)
= \\[1ex]
\nonumber
U_i^{n} - 0.5 \tau \left(V_{i j} \partial_x^{\kappa} U_{i j}^{n} + W_{i j} \partial_y^{\kappa} U_{i j}^{n} \right) + \\[1ex]
\label{ctusi2d1}
|\mathcal{C}_{ij} \mathcal{D}_{ij}|/12 \left( 2 U_{i j}^n + U_{i\pm1 j\pm1}^{n}+U_{i\mp1\, j\mp1}^{n}-U_{i+1\, j}^{n}-U_{i\, j+1}^n-U_{i-1\, j}^{n}-U_{i\, j-1}^{n} \right) ,
\end{eqnarray}
and
\begin{eqnarray}
\nonumber
U_{i j}^{n+1} + 
\tau V_{i j} \left(\partial_x^{\mp} U_{i j}^{n+1} - 0.5 \partial_x^{\kappa} U_{i\mp 1\, j}^{n+1} \right) + 
\tau W_{i j} \left(\partial_y^{\mp} U_{i j}^{n+1} - 0.5 \partial_y^{\kappa} U_{i\, j\mp 1}^{n+1} \right) + \\[1ex]
\nonumber
|\mathcal{C}_{ij} \mathcal{D}_{ij}|/6 \left( U_{i j}^{n+1}+U_{i\mp1\, j\mp1}^{n+1}-U_{i\mp1\, j}^{n+1}-U_{i\, j\mp1}^{n+1} \right)
= \\[1ex]
\nonumber
U_i^{n} - 0.5 \tau \left(V_{i j} \partial_x^{\kappa} U_{i j}^{n} + W_{i j} \partial_y^{\kappa} U_{i j}^{n} \right) - \\[1ex]
\label{ctusi2d2}
|\mathcal{C}_{ij} \mathcal{D}_{ij}|/12 \left( 2 U_{i j}^n + U_{i\mp1 j\pm1}^{n}+U_{i\pm1\, j\mp1}^{n}-U_{i+1\, j}^{n}-U_{i\, j+1}^n-U_{i-1\, j}^{n}-U_{i\, j-1}^{n} \right).
\end{eqnarray}
In (\ref{ctusi2d1}) and (\ref{ctusi2d2}) the identical convention is used for $\mp$ and $\pm$ as in (\ref{si2d}). 

Concerning the stability property, the numerical von Neumann stability analysis implies that the schemes (\ref{ctusi2d1}) and (\ref{ctusi2d2}) are unconditionally stable for the variable choice (\ref{kappasi2d}) of $\kappa$ parameters. This property is valid for the convex combination of these two schemes. 

As we show in the section on numerical experiments for the chosen representative examples, the scheme (\ref{ctusi2d1}) with (\ref{kappasi2d}) gives the most accurate results among all considered semi-implicit schemes. Moreover, the unconditional stability of this scheme is confirmed by these examples. 

\subsection{Implictly given computational domains}
\label{subsec-bc}
Our main motivation to introduce the unconditionally stable semi-implicit schemes is to solve the advection equation on computational domains with implicitly given boundaries. The basic idea is to use Cartesian grids even when the computational domain does not have a rectangular shape. In what follows we adopt an approach of an extrapolation of numerical solution for the grid nodes next to the implicitly given boundary $\partial \Omega$ as published in \cite{fmu14,fmu15}.

In particular, let $\phi=\phi(x,y)$, $(x,y) \in D$ be a given continuous function and $\Omega := \{(x,y) \in D, \phi(x,y) < 0\}$. We aim to solve the advection equation {\ref{adveq}) for $(x,y) \in \Omega \subset D$
using the Cartesian grid of the rectangular domain $D$ as described at the beginning of Section \ref{sec-2d}. We denote $\phi_{ij}:=\phi(x_i,y_j)$ for  $(x_i,y_j) \in D$ and search the values $U_{ij}^{n+1}$ only if $\phi_{ij}<0$.

Any scheme of the form (\ref{all2d}) can be used with no modifications if $\phi_{ij}<0$ and if for all its nonzero coefficients before $U_{i+k j}^{*}$ and $U_{i j+l}^{*}$ with $*=n, n+1$ one has that $\phi_{i+k j}\le 0$, resp. $\phi_{i j+l} \le 0$. Analogous considerations are valid for the CTU schemes (\ref{ctusi2d1}) or (\ref{ctusi2d2}) where one has to consider also the required diagonal values of numerical solution.

If the unmodified scheme (\ref{si2d}), (\ref{ctusi2d1}), or (\ref{ctusi2d2}) can not be used because some required neighbor values lie outside of $\Omega$, we exploit the advantage of variable choice for $\kappa$ parameters, and we use for the grid nodes next to the boundary $\partial \Omega$ only the scheme (\ref{si2d}) with $\kappa^x_{ij}=\hbox{sign}(\mathcal{C}_{ij})$ and $\kappa^y_{ij}=\hbox{sign}(\mathcal{D}_{ij})$. 
The motivation is that this scheme has the smallest stencil among all presented schemes, and it is unconditionally stable for any values of $\mathcal{C}_{ij}$ and $\mathcal{D}_{ij}$. The scheme  takes the form
\begin{eqnarray}
\label{si2dLW}
\nonumber
U_{ij}^{n+1} +  0.5 \, |\mathcal{C}_{ij}| ( U_{ij}^{n+1} - U_{i\mp 1 j}^{n+1} ) +  0.5 \, |\mathcal{D}_{ij}| ( U_{ij}^{n+1} - U_{i j\mp 1}^{n+1} ) = \\[1ex]
 U_{ij}^{n} -  0.5 \, |\mathcal{C}_{ij}| ( U_{i\pm 1 j}^{n} - U_{ij}^n )  -  0.5 \, |\mathcal{D}_{ij}| ( U_{i j\pm 1}^{n} - U_{ij}^n ) \,,
\end{eqnarray}
where the particular signs in $\mp$ or $\pm$ are chosen according to the signs of $\mathcal{C}_{ij}$ and $\mathcal{D}_{ij}$, see (\ref{si2d}).

The scheme (\ref{si2dLW}) can be used with no modifications if $\phi_{i\mp 1 j} \le 0$, $\phi_{i\pm 1 j} \le 0$, $\phi_{i j\mp 1} \le 0$, and $\phi_{i j\pm 1} \le 0$. We describe how to modify it if one has $\phi_{i\mp 1 j} > 0$ or $\phi_{i\pm 1 j} > 0$, the case $\phi_{i j\mp 1} > 0$ or $\phi_{i j\pm 1} > 0$ is treated analogously. 

Let $\phi_{i+1 j} > 0$, the case $\phi_{i-1 j} > 0$ is treated analogously. As $\phi_{ij} < 0$ one has that there exists a point $x_{i + \gamma} \in (x_i,x_{i+1})$, $\gamma \in (0,1)$ such that $x_{i + \gamma}=\gamma x_i + (1-\gamma) x_{i+1}$ and $\phi(x_{i+\gamma},y_j)=0$. One can determine the value $\gamma$ from an analytical form of $\phi$ or simply from the linear interpolation of $\phi_{ij}$ and $\phi_{i+1 j}$, see also \cite{fmu14,fmu15}. Note that in general the value of $\gamma$ can be arbitrary small and $x_{i+1} - x_{i+\gamma}=\gamma h$, so in a stability analysis the Courant number ${\cal C}_{ij}$ shall be divided  by $\gamma$.

If $\mathcal{C}_{ij}<0$ then the value $U_{i+1 j}^{n+1}$ is required by (\ref{si2dLW}). In this case the solution shall be given at the boundary node $(x_{i+\gamma},y_j)$, namely $u(x_{i+\gamma},y_j,t^{n+1})=u^D(x_{i+\gamma},y_j,t^{n+1})$, where the function $u^D$ is given. It corresponds to the case of inflow boundary with Dirichlet boundary conditions. Consequently one can replace the unavailable value $U_{i+1 j}^{n+1}$ in (\ref{si2dLW}) by the substitution
\begin{equation}
\label{subst}
U_{i+1 j}^{n+1} = \frac{1}{1-\gamma} \left( u^D(x_{i+\gamma},y_j,t^{n+1}) - \gamma U_{i j}^{n+1} \right)
\end{equation}
that is obtained simply by the linear extrapolation of the values $U_{i j}^{n+1}$ and $u^D(x_{i+\gamma},y_j,t^{n+1})$.

Furthermore, if $\mathcal{C}_{ij}>0$ then the value $U_{i+1 j}^{n}$ is required by (\ref{si2dLW}). This situation corresponds to the outflow boundary, when we apply the standard linear extrapolation \cite{lev02,fmu14,fmu15} to use the substitution
\begin{equation}
\label{outflow}
U_{i+1 j}^{n} = 2 U_{i j}^{n} -  U_{i-1 j}^{n} \,.
\end{equation}
In a rare situation when also the value $U_{i-1 j}^{n}$ is not available, we use the constant extrapolations $U_{i+1 j}^{n}=U_{i-1 j}^{n}=U_{ij}^n$.

\section{Numerical experiments}
\label{sec-num}
In what follows we study the properties of semi-implicit $\kappa$-schemes for some representative examples. In all examples the resulting linear algebraic systems are solved by Gauss-Seidel iterations using a strategy of so-called fast sweeping method where 
one sweep is given by four Gauss-Seidel iterations in four different orders  \cite{zao04,fmu15}. We use typically $1$ or $2$ sweeps that we specify for each example.
In all examples the domain is $D=(-1,1)^2$. The discretization steps are $h=2/M$ and $\tau=T/N$ where $M$, $T$ and $N$ will be given. The maximal Courant number is defined as the maximum among all directional grid Courant numbers $|{\cal C}_{ij}|$ and $|{\cal D}_{ij}|$ given in (\ref{cn2d}).

To check the implementation we computed examples with $u^0(x,y)$ being a randomly chosen quadratic function and $\vec{V}$ being an arbitrary constant velocity. For all considered numerical schemes the exact solution is reproduced by the numerical solution up to a machine accuracy for any chosen $h$ and $\tau$. Choosing the initial function as a cubic polynomial, only the convex combination of CTU schemes (\ref{ctusi2d1}) and (\ref{ctusi2d2}) with the variable choice (\ref{kappasi2d}) of $\kappa$ parameters 
gives numerical solutions differing from the exact ones purely by rounding errors.

\subsection{Implicitly given computational domain}

To illustrate the advantages of semi-implicit schemes for the applications on implicitly given domains, we solve the advection equation (\ref{adveq}) only inside of a circle of radius $1$ given implicitly as a zero level set of the function $\phi(x,y)=\sqrt{x^2+y^2}-1$ for $(x,y) \in D$ following the approach of subsection \ref{subsec-bc}.
We compute the examples with the advection velocity defined by 
\begin{equation}
\label{rotation}
\vec{V}=(-2 \pi y \, , \, 2 \pi x ) \,.
\end{equation}
The exact solution for any initial function $u^0$ is given by $u(x,y,t)=u^0(x \cos (2 \pi t) + y \sin(2 \pi t) , y \cos(2 \pi t) - x \sin (2 \pi t))$. We consider $T=1$, so the initial profiles given by $u^0$ rotate once to return at $t=1$ to their initial position. We consider two typical initial level set functions - the first one representing implicitly a smooth interface  (e.g. a circle \cite{Min2007,fm07,Phongthanapanich2010,Kim2011,Wang2012,Samuel2014, Starinshak2014}), and the second one representing a piecewise smooth interface  \cite{for05,Min2007,fm07,Phongthanapanich2010,Starinshak2014} (e.g. a square). 

To estimate the EOCs we compute the error
\begin{equation}
\label{error}
E = h^2 \max_{n=1,..,N} \sum_{i,j=1}^M  |U_{i j }^n - u(x_i,y_j,t^n)| \,.
\end{equation}
We choose in all examples the Cartesian grids with $M=40,80,160$. The minimal values of $\gamma$ in (\ref{subst}) for these three Cartesian grids are roughly $0.079$, $0.044$, and $0.025$, so a certain number of grid points next to the circular boundary has a neighbor point at the boundary in a distance approximately $\gamma h$. Consequently, the local grid Courant numbers ${\cal C}_{ij}$ or ${\cal D}_{ij}$ corresponding to those grid points are multiplied by 
the values $1/\gamma$, so they become effectively between $12$ and $40$ times larger. 

The results are compared for the following choices of $\kappa$ parameters in the scheme (\ref{si2d}) without CTU extension, see also (\ref{siLW}), (\ref{siF}), and (\ref{siQ}),
\begin{eqnarray}
\label{kp}
\kappa^x = \hbox{sign}(\mathcal{C}_{ij}) \,, \quad \kappa^y = \hbox{sign}(\mathcal{D}_{ij}) \,, \\
\label{km}
\kappa^x = -\hbox{sign}(\mathcal{C}_{ij}) \,, \quad \kappa^y = - \hbox{sign}(\mathcal{D}_{ij}) \,, \\
\label{k0}
\kappa^x = \kappa^y = 0 \,,\\
\label{k3}
\kappa^x = \hbox{sign}(\mathcal{C}_{ij})(1-|\mathcal{C}_{ij})/3 \,, \quad \kappa^y = \hbox{sign}(\mathcal{D}_{ij})(1-|\mathcal{D}_{ij})/3 \,. 
\end{eqnarray}
The variable choice (\ref{k3}) is used additionally also with the CTU scheme (\ref{ctusi2d1}).

We begin with the initial function $u^0(x,y)=\sqrt{x^2+(y-0.5)^2}$ being the distance function in Euclidian metric. It  can be viewed as an example of the level set function representing a circular interface \cite{Min2007,fm07,Phongthanapanich2010,Kim2011,Wang2012,Samuel2014, Starinshak2014} (a smooth one). 

Firstly, we compute the example with (\ref{si2d}) and (\ref{kp} - (\ref{k3}) using $N=5 M/4$ (the maximal Courant number far from the circular boundary being $2.5$). 
The results are summarized in Table \ref{tab02} and in Figure \ref{fig02a}.  The schemes show the EOC approaching the $2^{nd}$ order accuracy from below in this example. No instabilities in the numerical solutions are observed for the cut cells. To solve the linear systems one sweep was enough.

\begin{figure}[h!]
\begin{center}
\includegraphics[width=0.24\columnwidth]{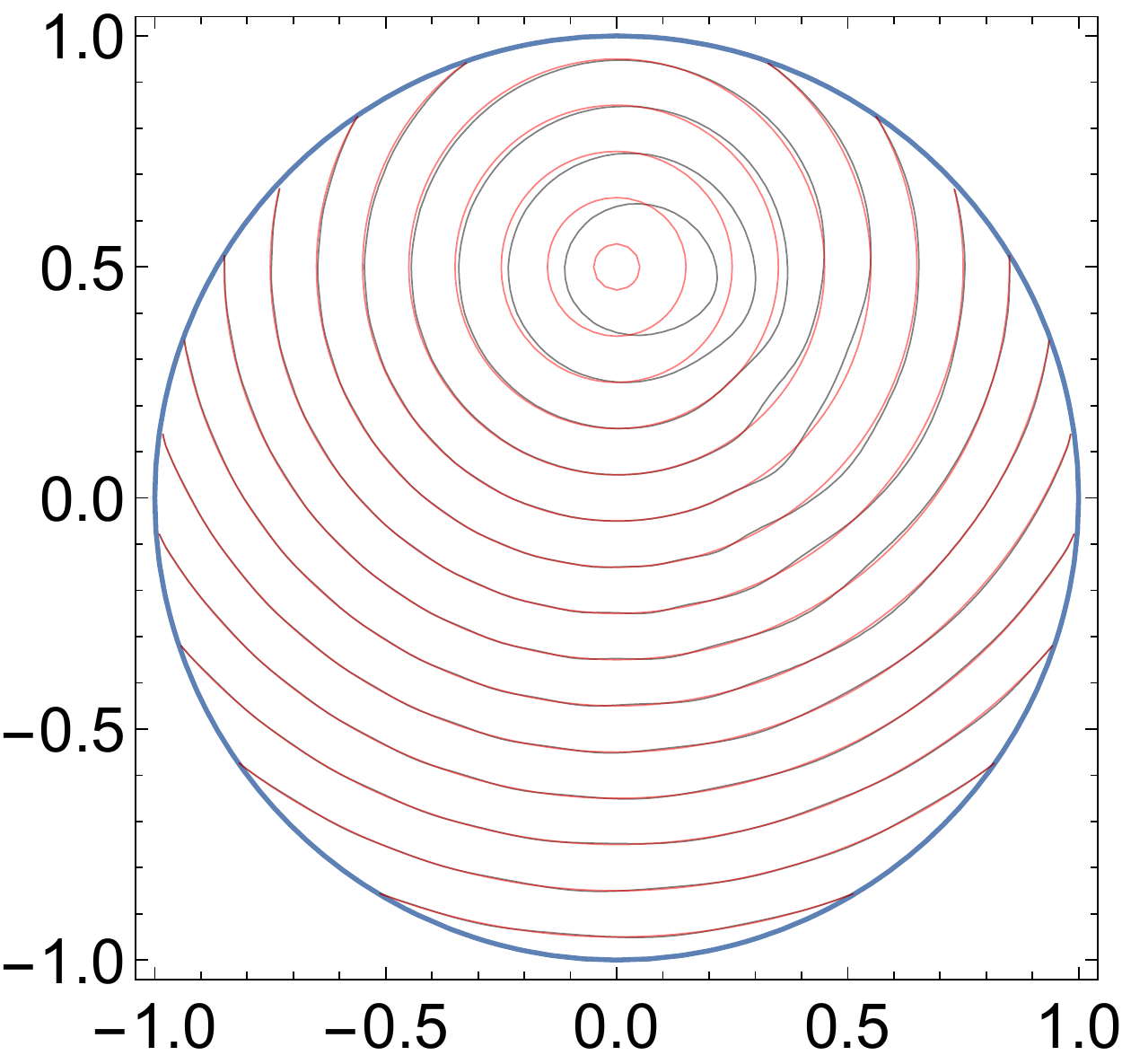}
\includegraphics[width=0.24\columnwidth]{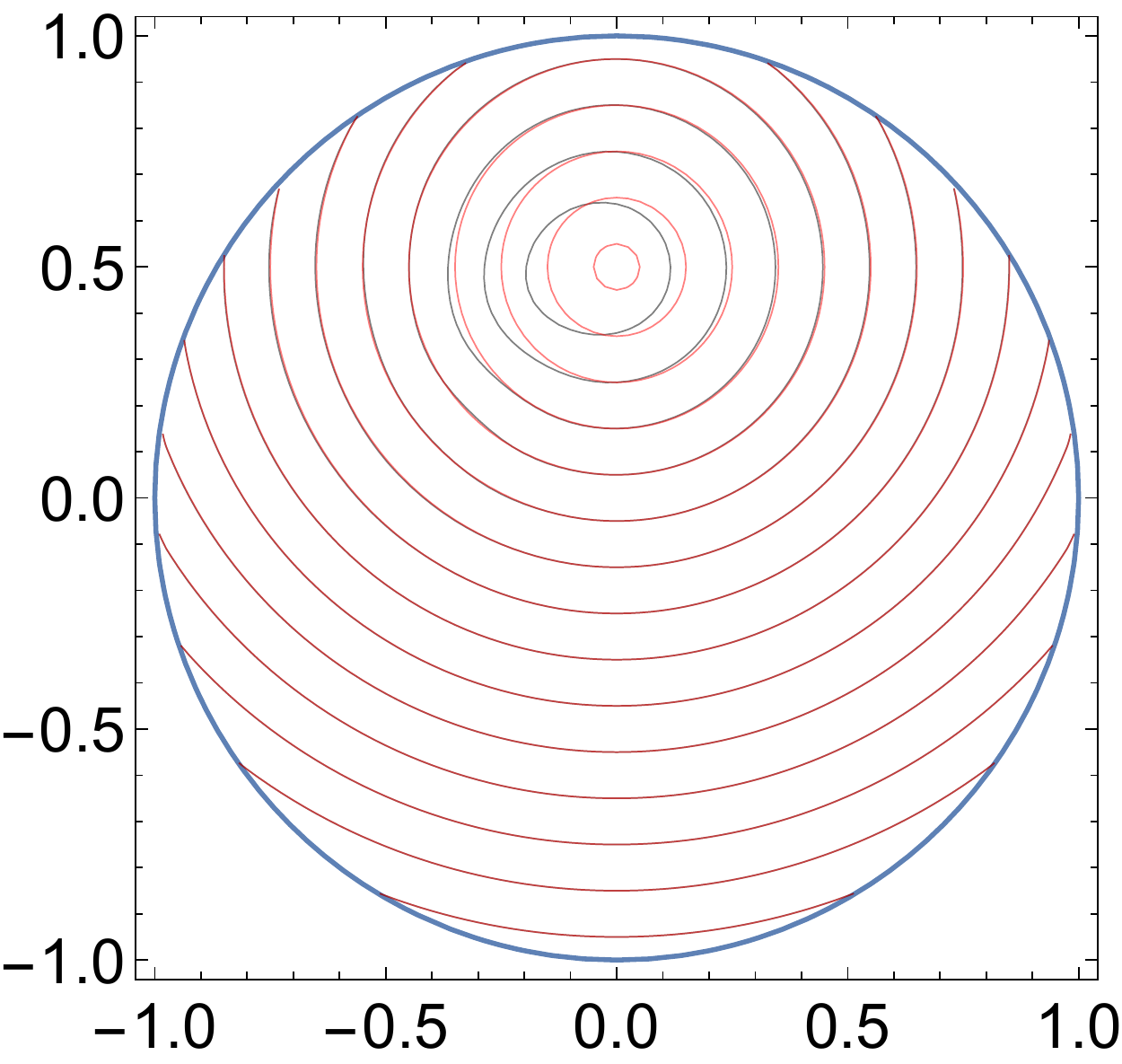}
\includegraphics[width=0.24\columnwidth]{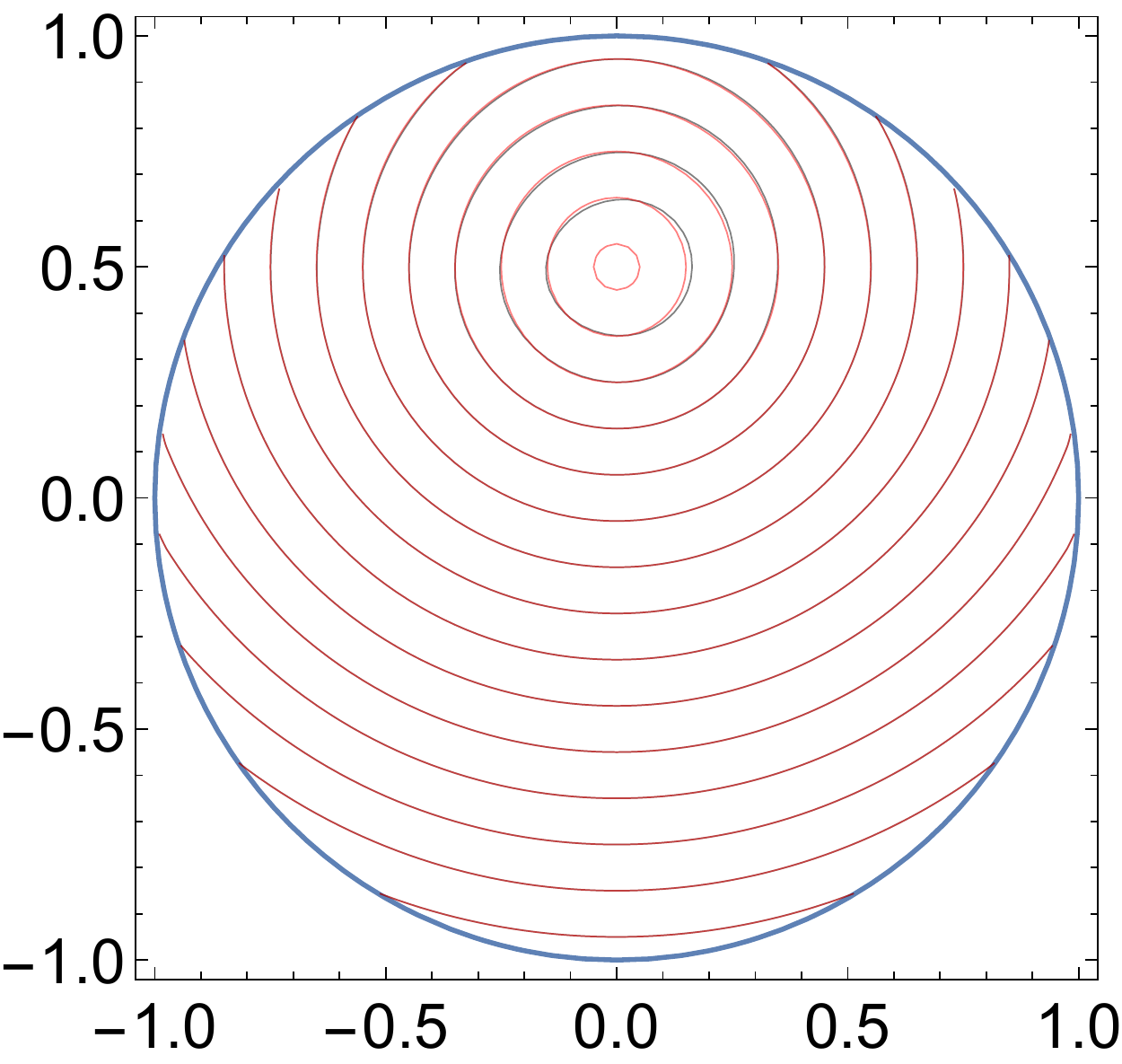}
\includegraphics[width=0.24\columnwidth]{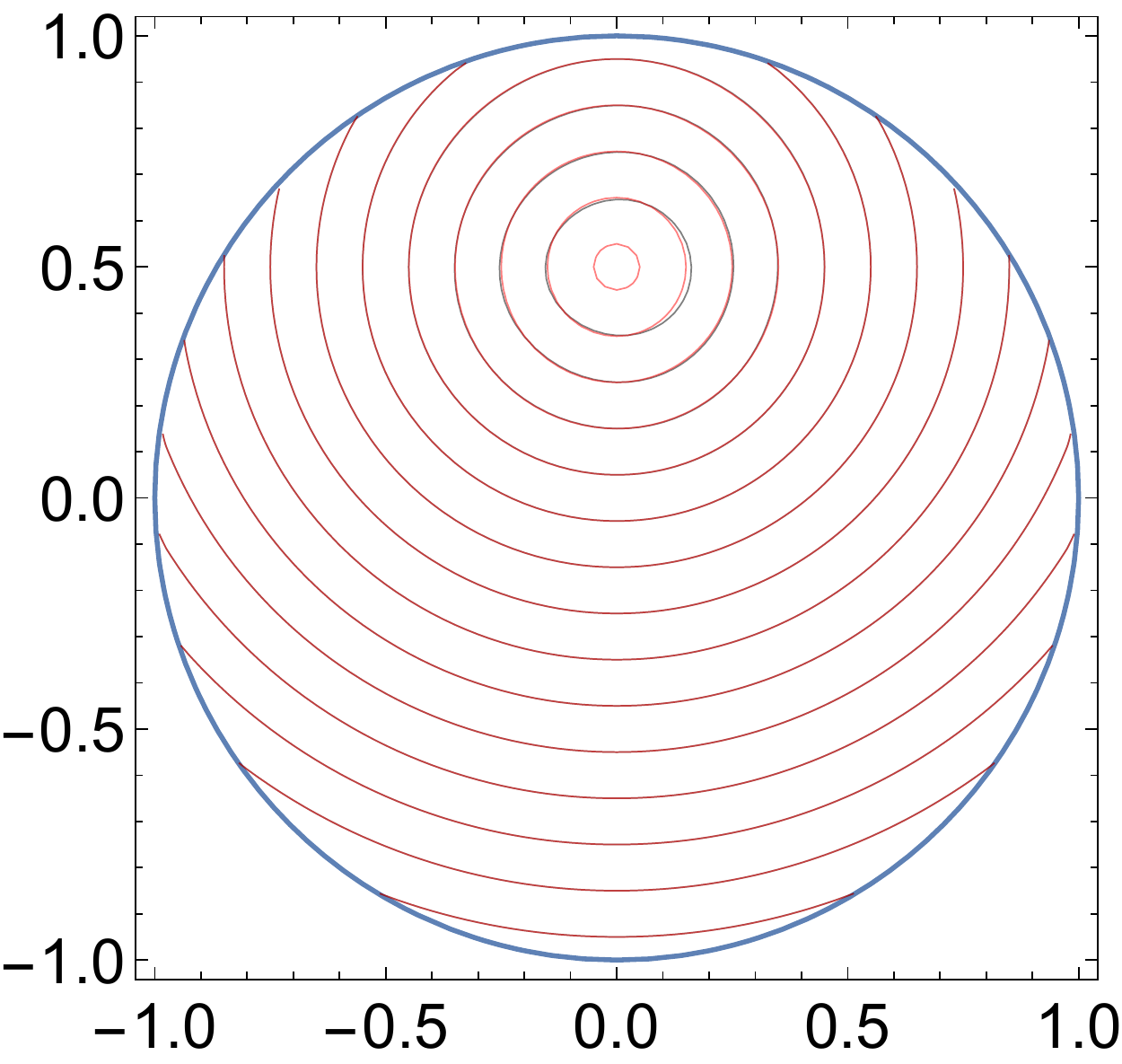}
\caption{Numerical solutions for the rotation of distance function in Euclidian metric for $M=80$ at $t=1$  using (\ref{si2d}) with (\ref{kp})  (the first picture), (\ref{km}), (\ref{k0}), and (\ref{k3}) (the last picture). The chosen time step corresponds to the maximal Courant number for the grid nodes far from the circular boundary to $2.5$. The red contour lines represent the exact solution, the black contour lines represent the numerical solutions for values $0.1$, $0.2$ up to $1.4$. Note a shift in the position of the contour line $0.2$ in the numerical solutions of schemes (\ref{kp}) and (\ref{km}) that is called a ``phase error'' \cite{lev02} in the case of analogous fully explicit schemes .
}
\label{fig02a}
\end{center}
\end{figure}

\begin{table}[h!]
\begin{center}
\begin{tabular}{||c||c|c||c|c||c|c||c|c||}
\hline
$M$ & $E_{ci}$ & $E_{sq}$ & $E_{ci}$ & $E_{sq}$ & $E_{ci}$ & $E_{sq}$ & $E_{ci}$ & $E_{sq}$  \\
\hline
40  & 47.2 & 120. & 27.3 & 66.4 & 9.47 & 47.7 & 6.54 & 35.1 \\
80  & 15.1 & 62.3 & 8.34 & 32.8 & 2.78 & 23.2 & 1.77 & 15.0\\
160& 4.58 & 31.4 & 2.46 & 15.2 & .782 & 9.92 & .484 & 5.82 \\
\hline
\end{tabular}
\caption{The error (\ref{error}) (multiplied by $10^3$) for the rotation of distance function in Euclidian metric (the columns with $E_{ci}$) and the rotation of distance function in maximum metric (the columns with $E_{sq}$). The examples are computed in the implicitly given circular domain using (\ref{si2d}) with (\ref{kp}) (the first two columns), (\ref{km}), (\ref{k0}), and (\ref{k3}) (the last two columns) with $N=5 M/4$, (the maximal Courant number far from the circular boundary being $2.5$). Note that the schemes are approaching the $2^{nd}$ order accuracy from below for the smooth case and the $1^{st}$ order accuracy from above for the non-smooth case.}
\label{tab02}
\end{center}
\end{table}

One can observe a so-called ``phase error'' \cite{lev02} for this example in the form of a shift for the location of some contour lines with respect to the exact position when using the schemes (\ref{kp}}) and (\ref{km}). These schemes are based  on the one sided approximations of gradient in (\ref{grad2d}) - (\ref{grad2dB}), whereas the scheme (\ref{k0}) uses the central approximation for which the phase error is much less visible. The scheme $(\ref{k3})$ gives the most accurate results for this example. We remind that we use the scheme (\ref{kp}) locally as described in subsection \ref{subsec-bc}  for the grid nodes next to the circular boundary. 

\begin{figure}[h!]
\begin{center}
\includegraphics[width=0.24\columnwidth]{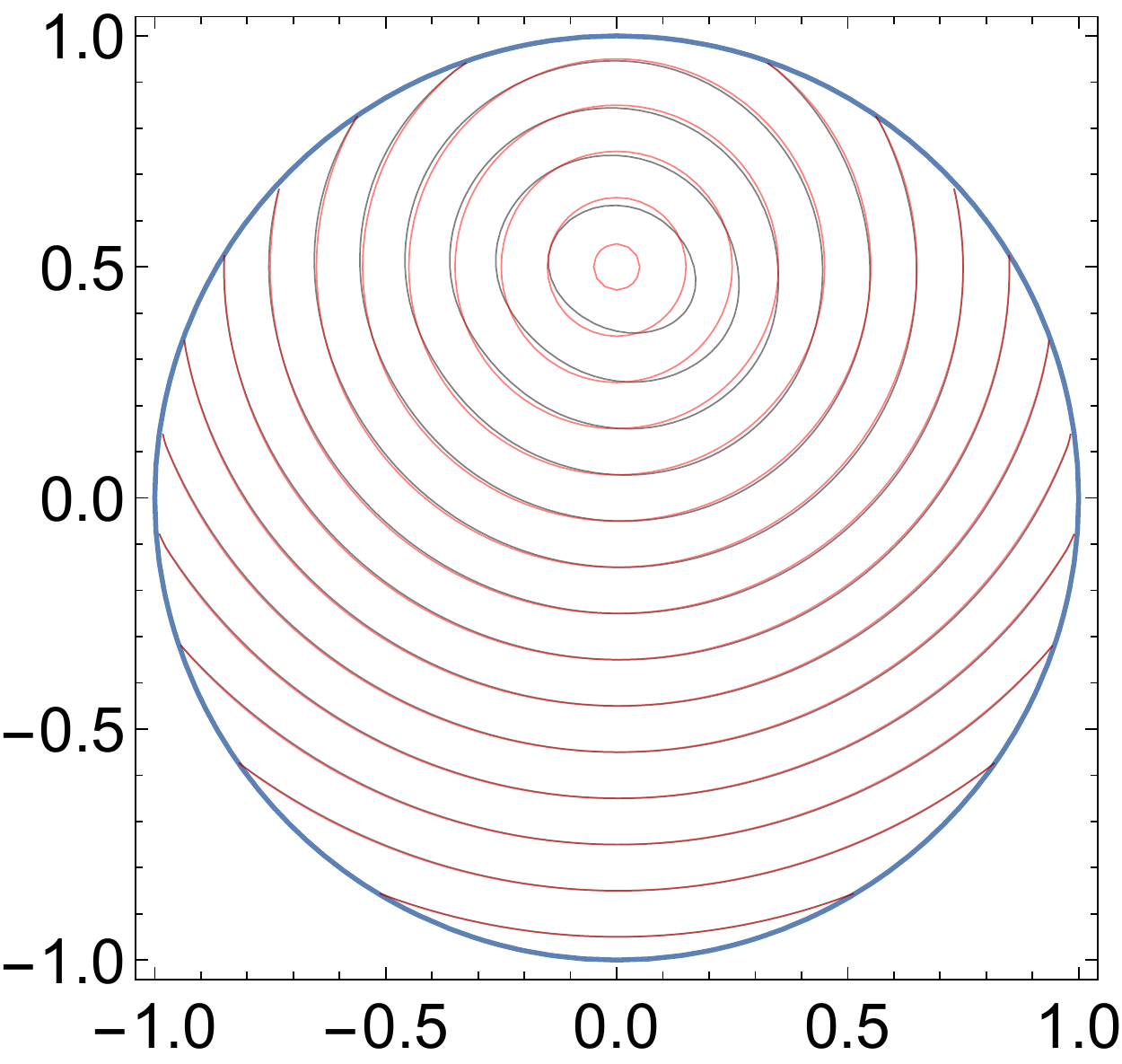}
\includegraphics[width=0.24\columnwidth]{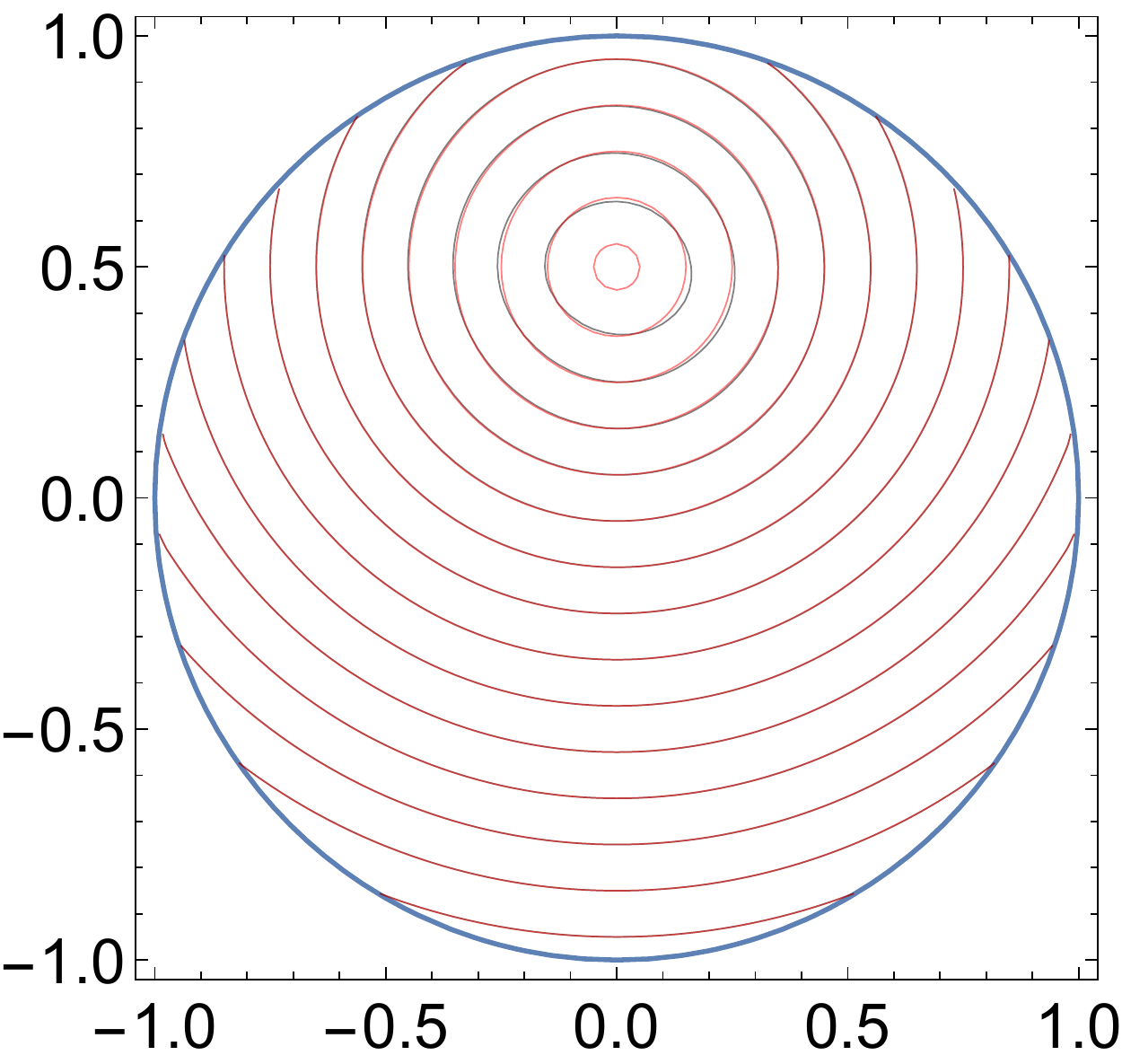}
\includegraphics[width=0.24\columnwidth]{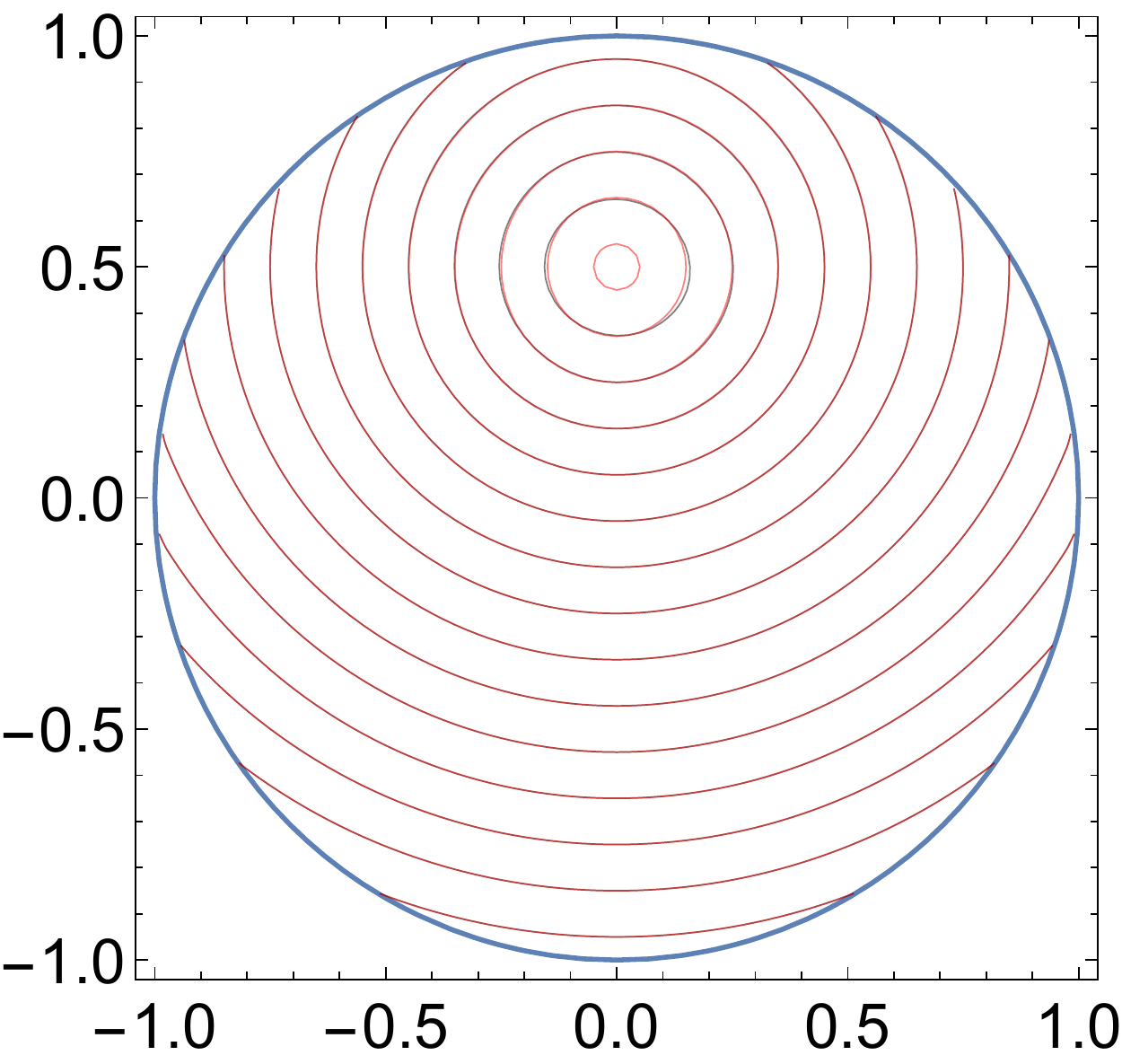}
\includegraphics[width=0.24\columnwidth]{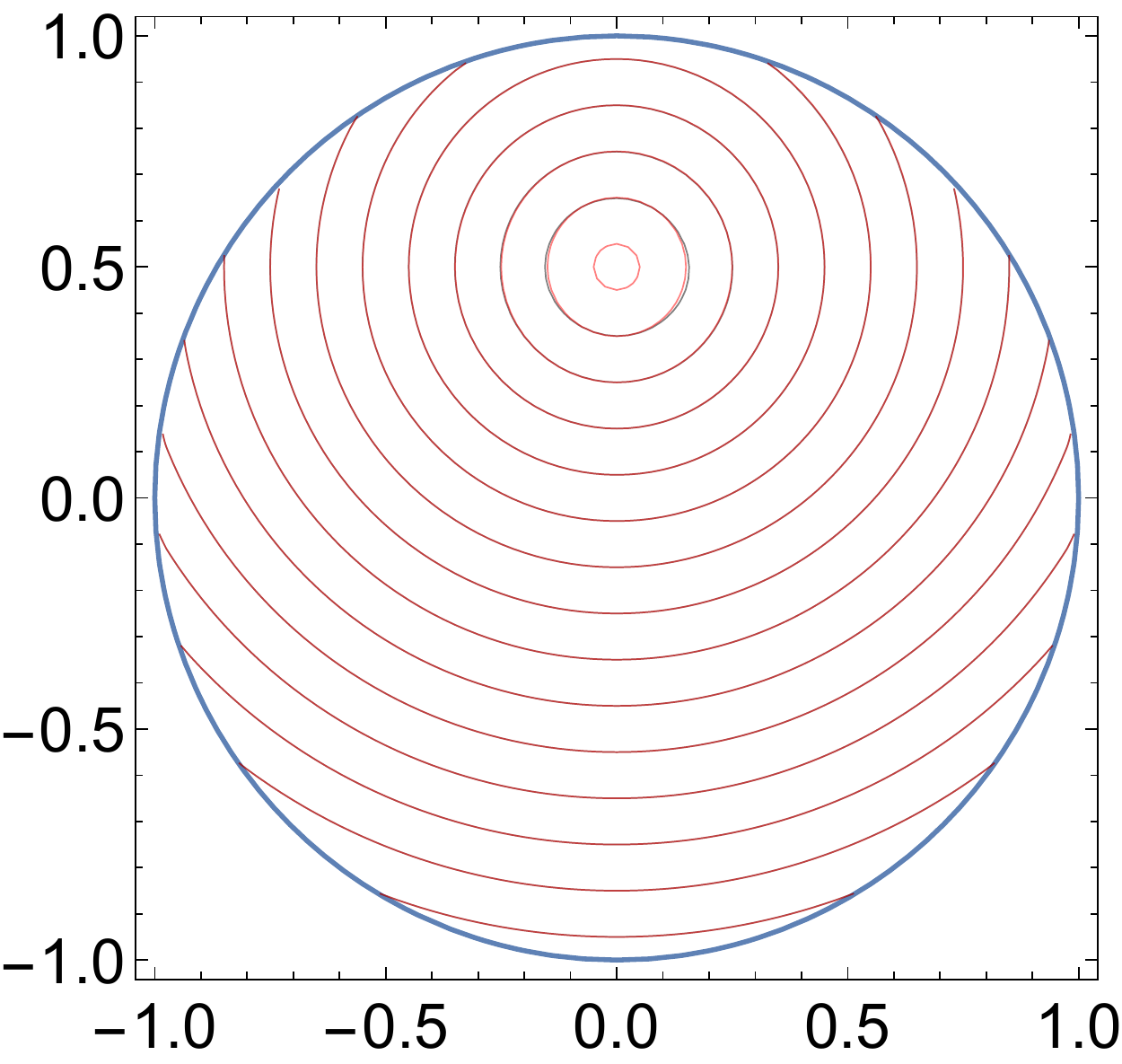}
\caption{Numerical solutions for the rotation of distance function in Euclidian metric in the implicitly given circular domain for $M=80$ at $t=1$  using (\ref{ctusi2d1} and (\ref{k3}) with the time steps corresponding to the maximal Courant numbers for the grid nodes far from the circular boundary to $10$ (the first picture), $5$, $2.5$ and $1.25$ (the last picture). The red contour lines represent the exact solution, the black contour lines represent the numerical solutions for the values $0.1$, $0.2$ up to $1.4$.
}
\label{fig02b}
\end{center}
\end{figure}

Furthermore, we compute the example on the medium grid with $M=80$ using the CTU scheme (\ref{ctusi2d1}) with (\ref{k3}) for $N=25, 50, 100, 200$ that corresponds to the maximal Courant numbers far from the circular boundary having the values $10, 5, 2.5, 1.25$. The results are presented in Figure \ref{fig02b}.
The error (\ref{error}) takes the values $12.6$, $4.20$, $1.74$, and $1.04$ multiplied by $10^{-3}$. To solve the linear systems one sweep was used except for the case $N=25$ when two sweeps were used.

For the next example we choose $u^0(x,y)= \max\{|x+0.5|,|y|\}$, i.e. the distance function in the maximum metric. It can be viewed as an example of the level set function representing a squared interface \cite{for05,Min2007,fm07,Phongthanapanich2010,Starinshak2014} (i.e. a piecewise smooth interface in this case). We compute the example with (\ref{si2d}) and (\ref{kp}) - {\ref{k3}) using $N=5 M/4$, see Table \ref{tab02}, and $N=5 M/2$, see Figure \ref{fig03a}, i.e. the maximal Courant number for the grid nodes far from the circular boundary being $2.5$ and $1.25$, respectively.
The schemes show the EOC approaching the $1^{st}$ order accuracy from above in this example, i.e  for the non-smooth solution. No instabilities in the numerical solutions are observed for the cut cells and one sweep was enough to solve the linear systems.

\begin{figure}[h!]
\begin{center}
\includegraphics[width=0.24\columnwidth]{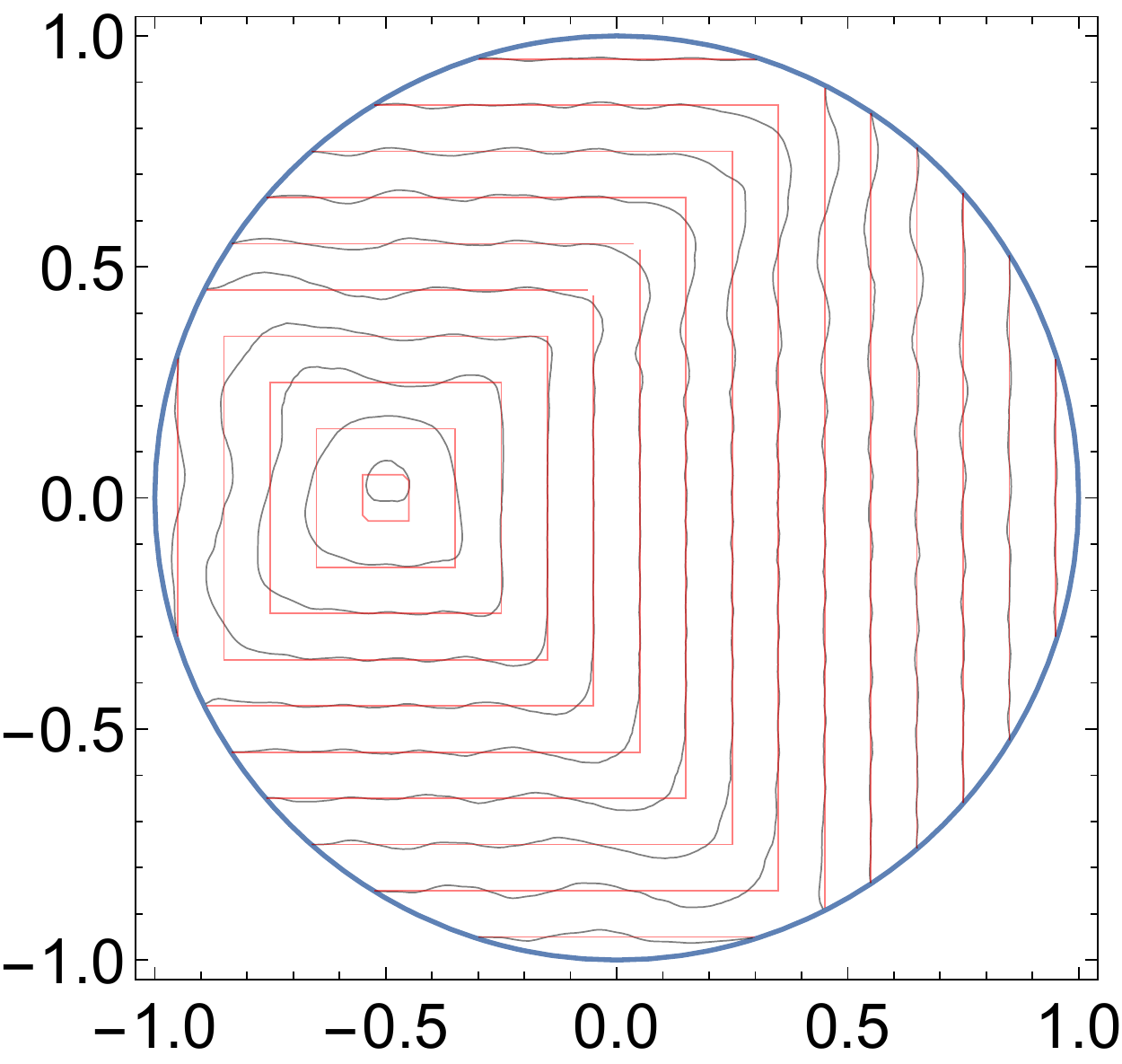}
\includegraphics[width=0.24\columnwidth]{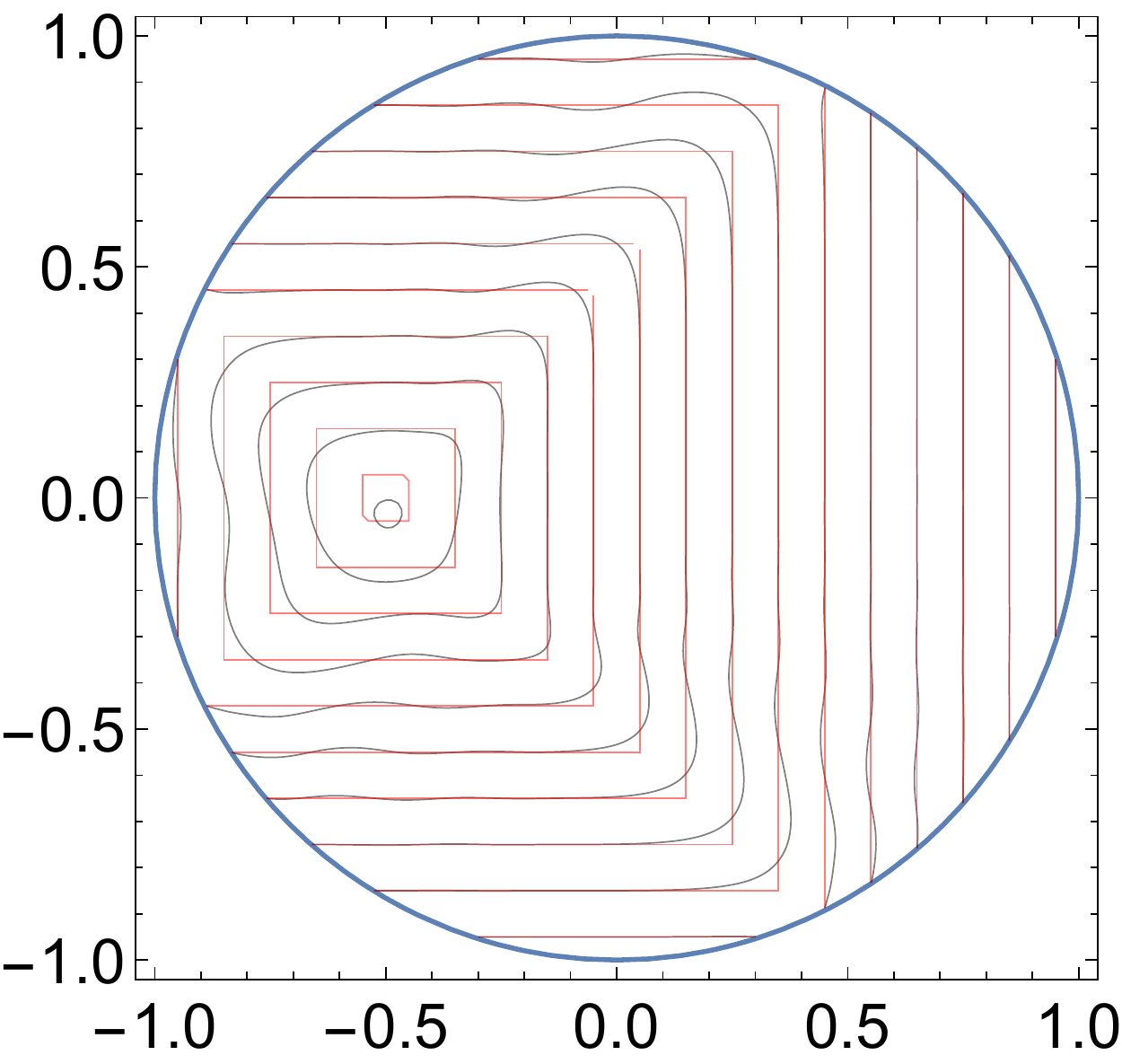}
\includegraphics[width=0.24\columnwidth]{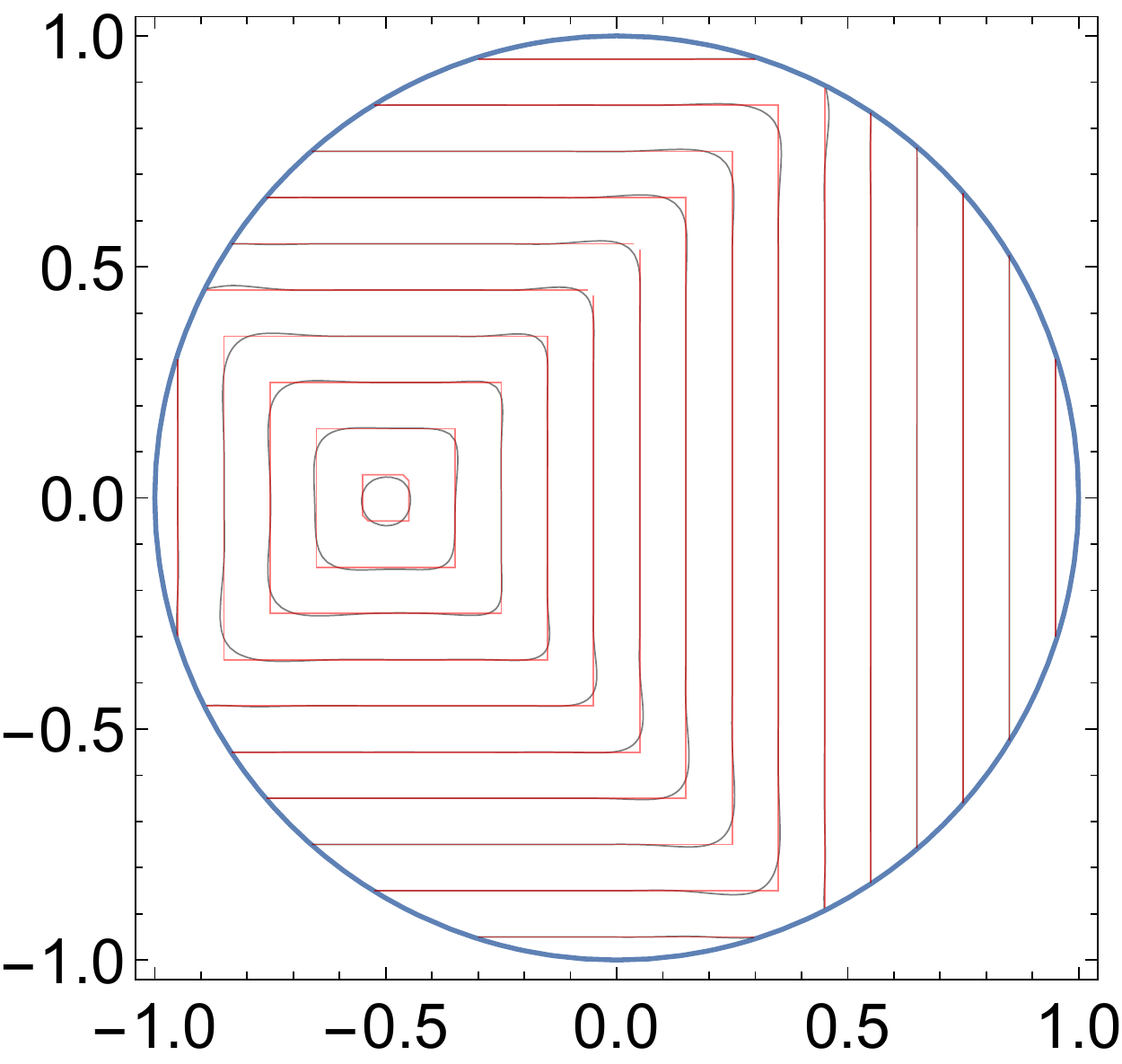}
\includegraphics[width=0.24\columnwidth]{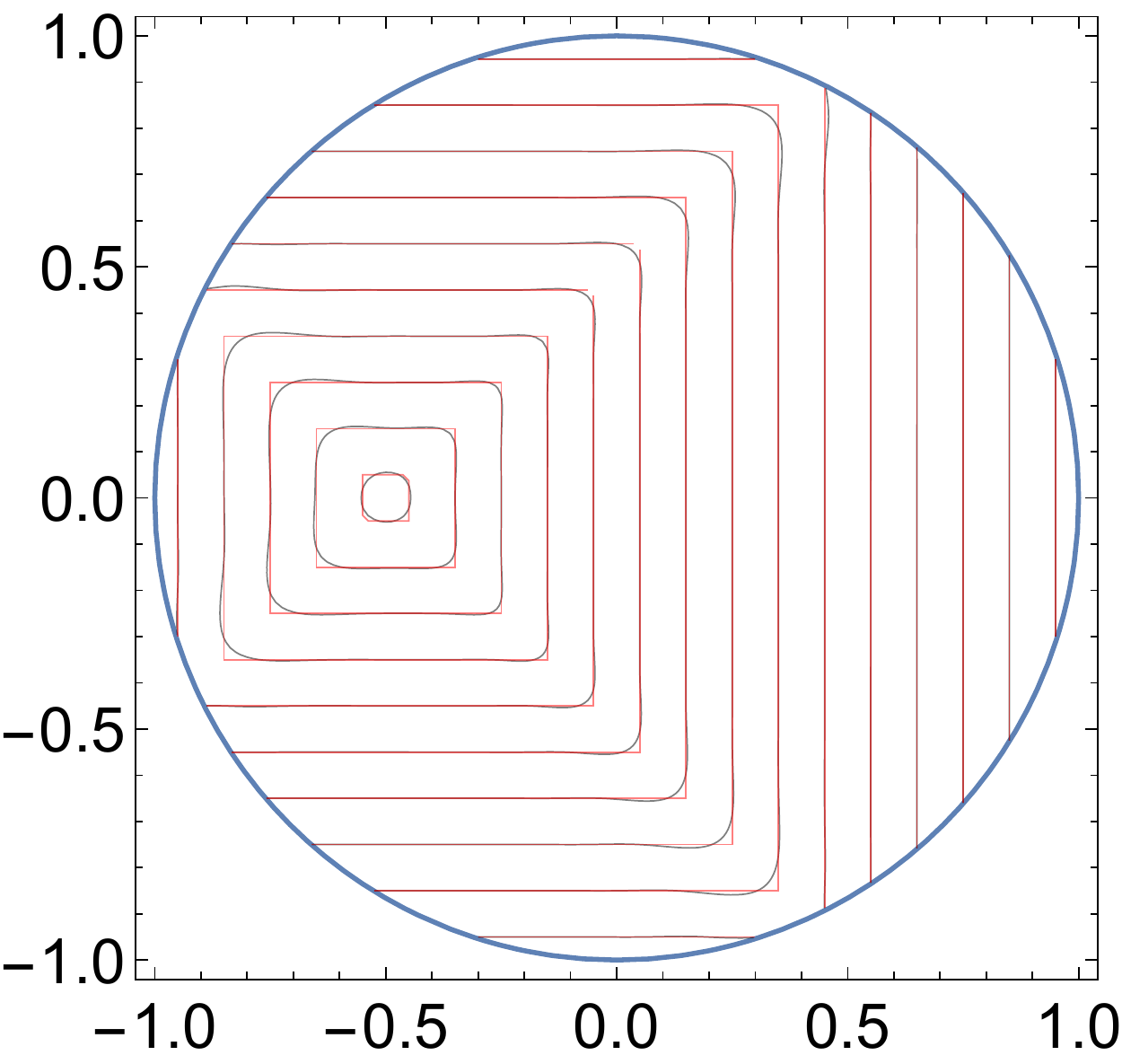}
\caption{Numerical solutions for the rotation of a distance function in the maximum metric for $M=160$ at $t=1$  using (\ref{si2d}) with (\ref{kp})  (the first picture), (\ref{km}), (\ref{k0}), and (\ref{k3}) (the last picture). The chosen time step corresponds to the maximal Courant number for the grid nodes far from the circular boundary to $1.25$. The red contour lines represent the exact solution, the black contour lines represent the numerical solutions for the values $0.1$, $0.2$ up to $1.4$. Note that the schemes (\ref{kp}) and (\ref{km}) based on the one sided finite difference approximations of gradient are visibly less accurate than the scheme (\ref{k0}) based on the central difference approximation. Moreover, the oscillations in numerical solution for (\ref{kp}) are not damped as the amplification factor for this schemes is $1$.
}
\label{fig03a}
\end{center}
\end{figure}

One can again observe the phase error \cite{lev02} for this example when using the schemes (\ref{kp}}) and (\ref{km}) that are based  on the one sided approximations of gradient in (\ref{grad2d}) - (\ref{grad2dB}). Moreover some  oscillations of contour lines for the numerical solution obtained with (\ref{kp}) are not damped as the amplification factor $|S|$ for this scheme is equal $1$ everywhere. The scheme $(\ref{k3})$ gives the most accurate results for this example. 

\begin{figure}[h!]
\begin{center}
\includegraphics[width=0.24\columnwidth]{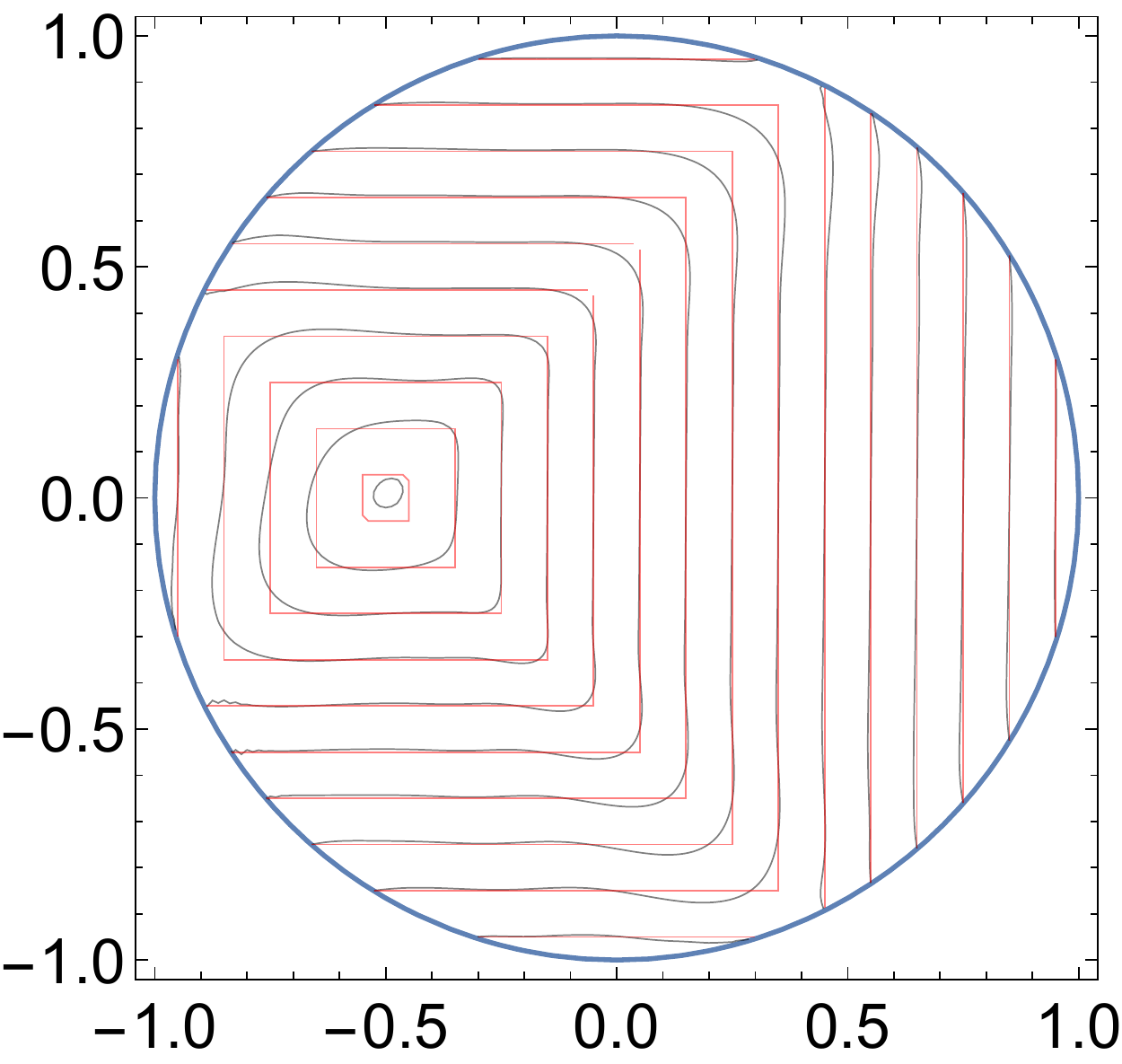}
\includegraphics[width=0.24\columnwidth]{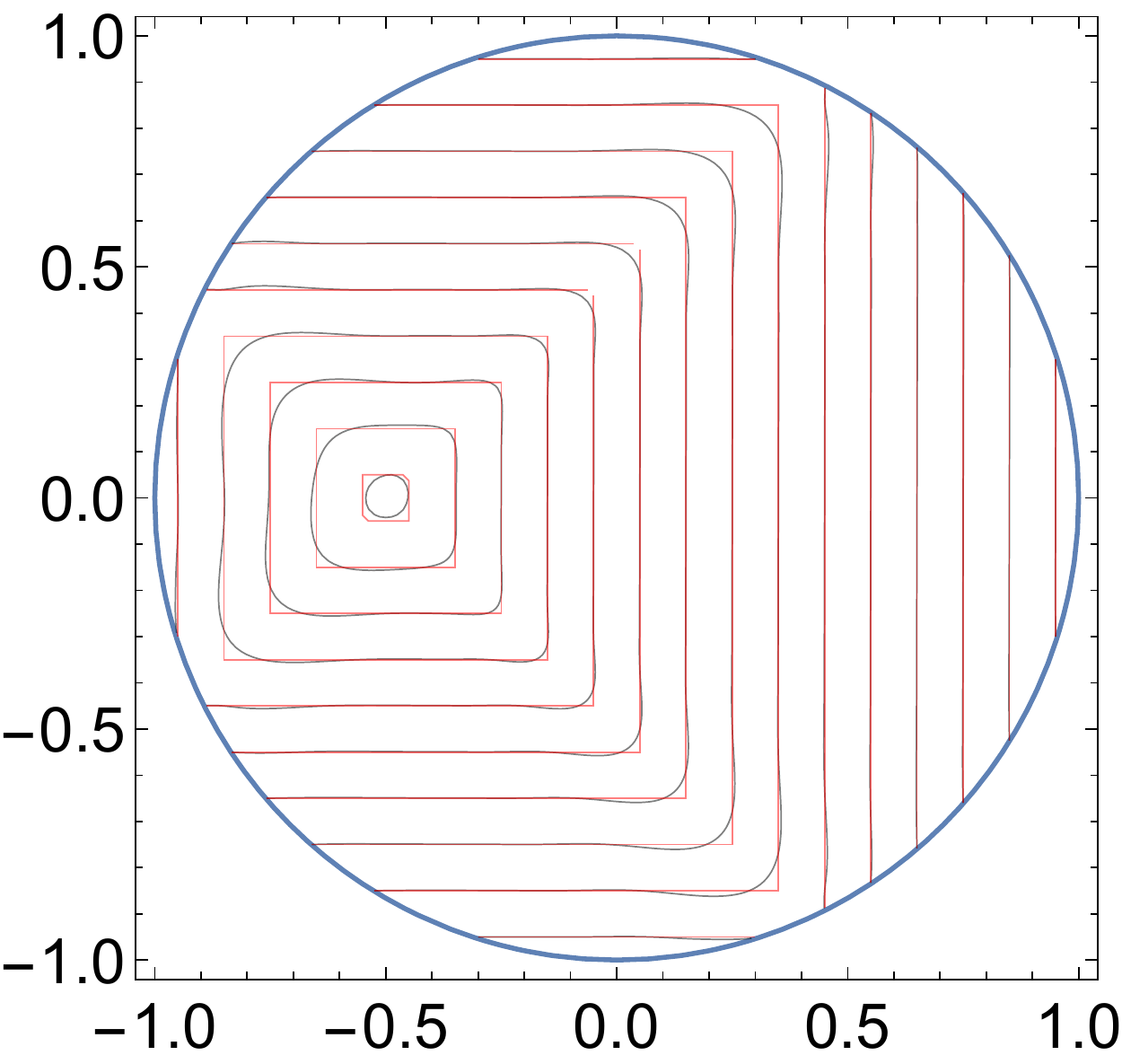}
\includegraphics[width=0.24\columnwidth]{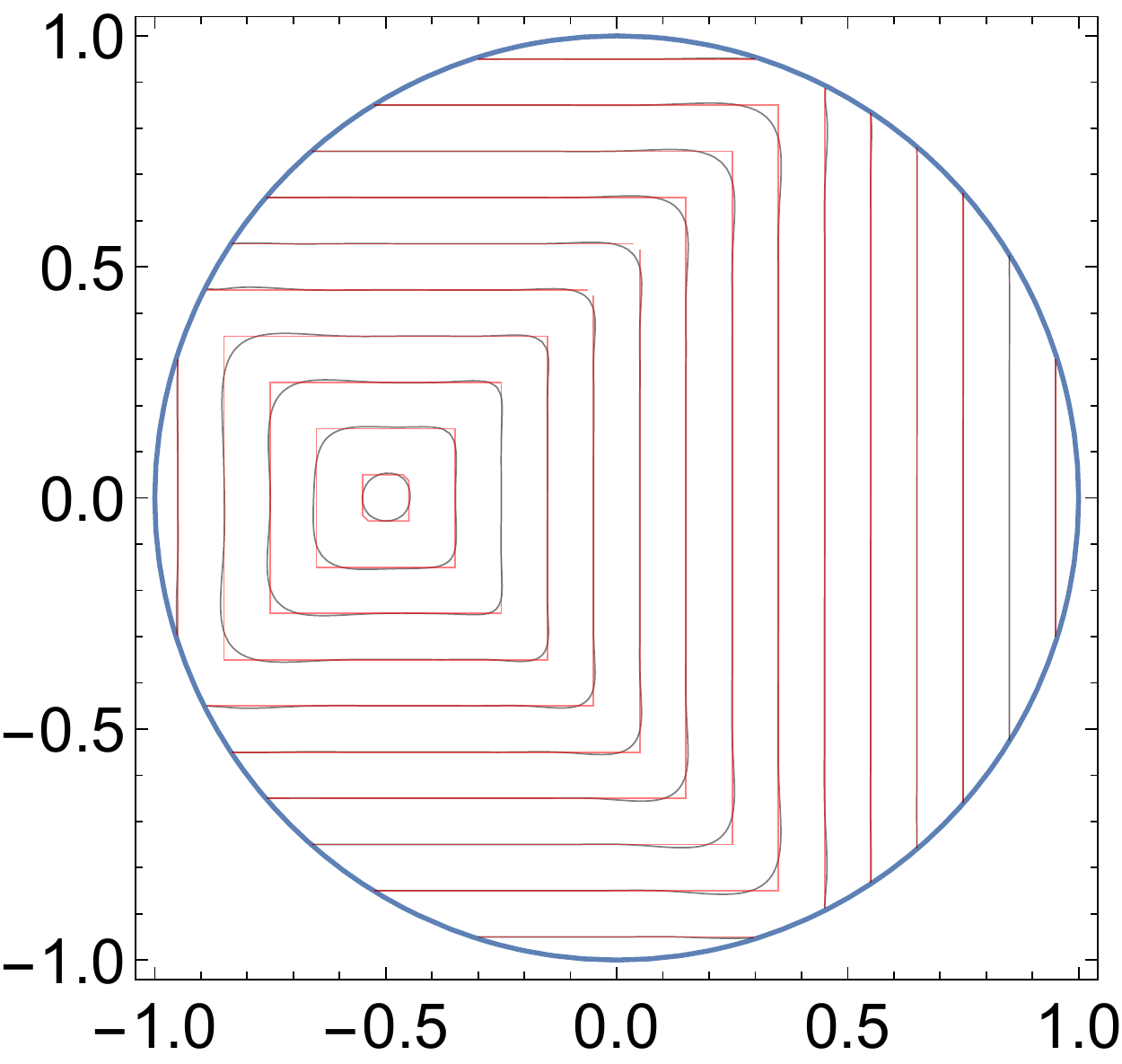}
\includegraphics[width=0.24\columnwidth]{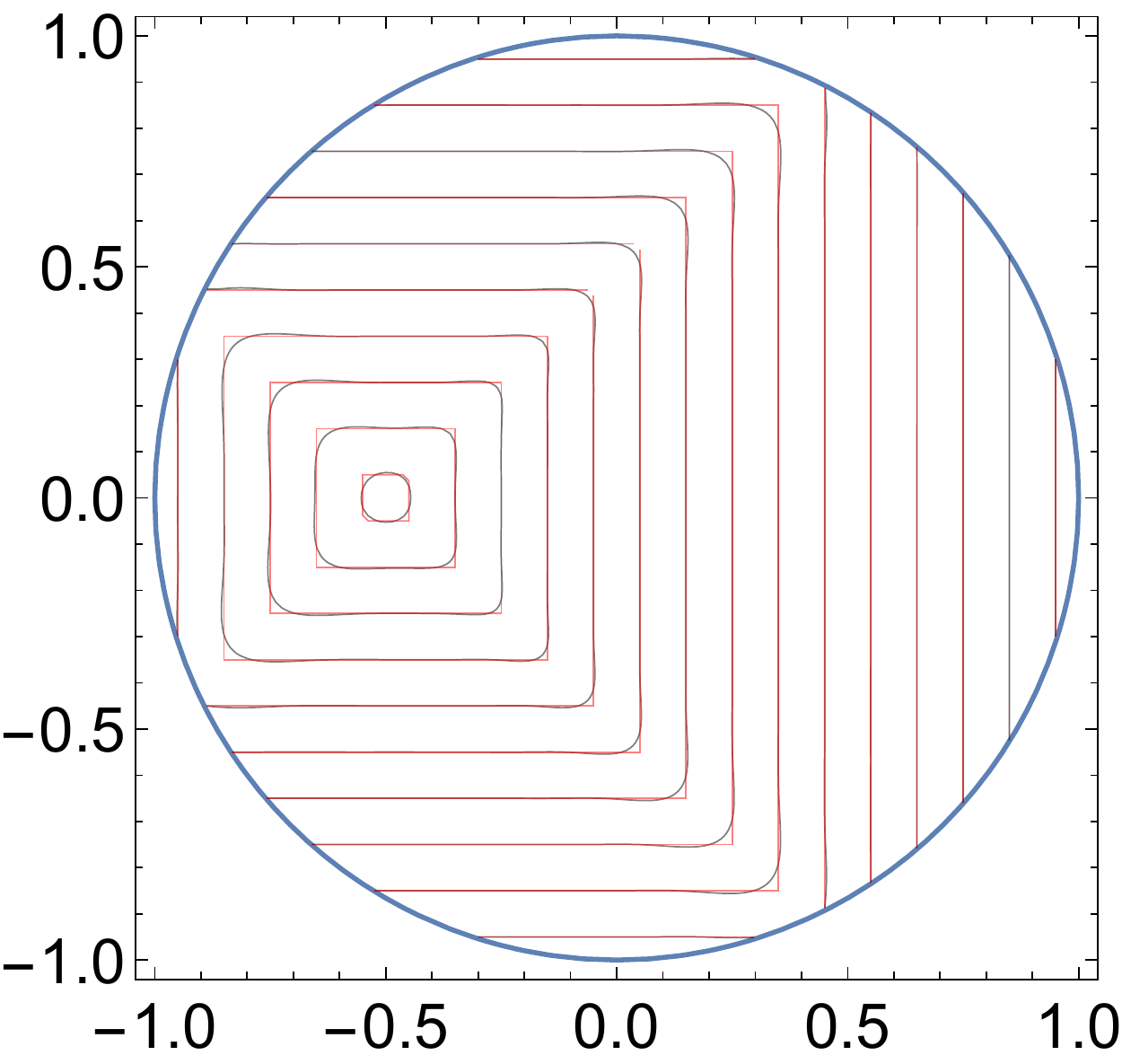}
\caption{Numerical solutions for the rotation of the distance function in the maximum metric in the implicitly given circular domain for $M=160$ at $t=1$  using (\ref{ctusi2d1} and (\ref{k3}) with the time steps corresponding to the maximal Courant numbers for the grid nodes far from the circular boundary to $10$ (the first picture), $5$, $2.5$ and $1.25$ (the last picture). The red contour lines represent the exact solution, the black contour lines represent the numerical solutions for values $0.1$, $0.2$ up to $1.4$.
}
\label{fig03b}
\end{center}
\end{figure}

Finally, we compute the example on the finest grid with $M=160$ using the CTU scheme (\ref{ctusi2d1}) with (\ref{k3}) for $N=50, 100, 200, 400$ that corresponds to the maximal Courant numbers far from the circular boundary having the values $10, 5, 2.5, 1.25$. The results are presented in Figure \ref{fig03b} .
The error (\ref{error}) takes the values $18.7$, $8.01$, $4.86$, and $3.54$ multiplied by $10^{-3}$. The linear systems were solved using only one sweep.

\subsection{Example with largely varying velocity}

In the next example we illustrate another advantage of semi-implicit schemes for the solution of advection equation when the velocity is varying significantly in the computational domain. 
To do so we choose the velocity varying exponentially in a diagonal direction,
\begin{equation}
\label{vexp}
\vec{V} = (e^{2(y-x)},e^{2(y-x)}) \,.
\end{equation}
The initial function $u^0(x,y)$ is given by the Euclidian distance to $(-1,-1)$. We fix the values $u(x,y,t)=u^0(x,y)$ for the inflow part of $D$, namely the east and the south sides of the squared domain $D$, the other two sides are outflow boundaries. The exact solution is given by the method of characteristics where the fixed values at inflow boundaries are respected. Note that the solution has discontinuous first derivatives with respect to $x$ an $y$. 

The speed $|\vec{V}|$ in the corner $(-1,1)$ is about $2981$ times larger than in the corner $(1,-1)$, so a large variation of $\vec{V}$ occurs in the domain $D$. The stationary solution for this example is the function $|y-x|$ and these values are approached by the time dependent solution $u(x,y,t)$ very rapidly in the left top corner of the domain.

\begin{figure}[h!]
\begin{center}
\includegraphics[width=0.24\columnwidth]{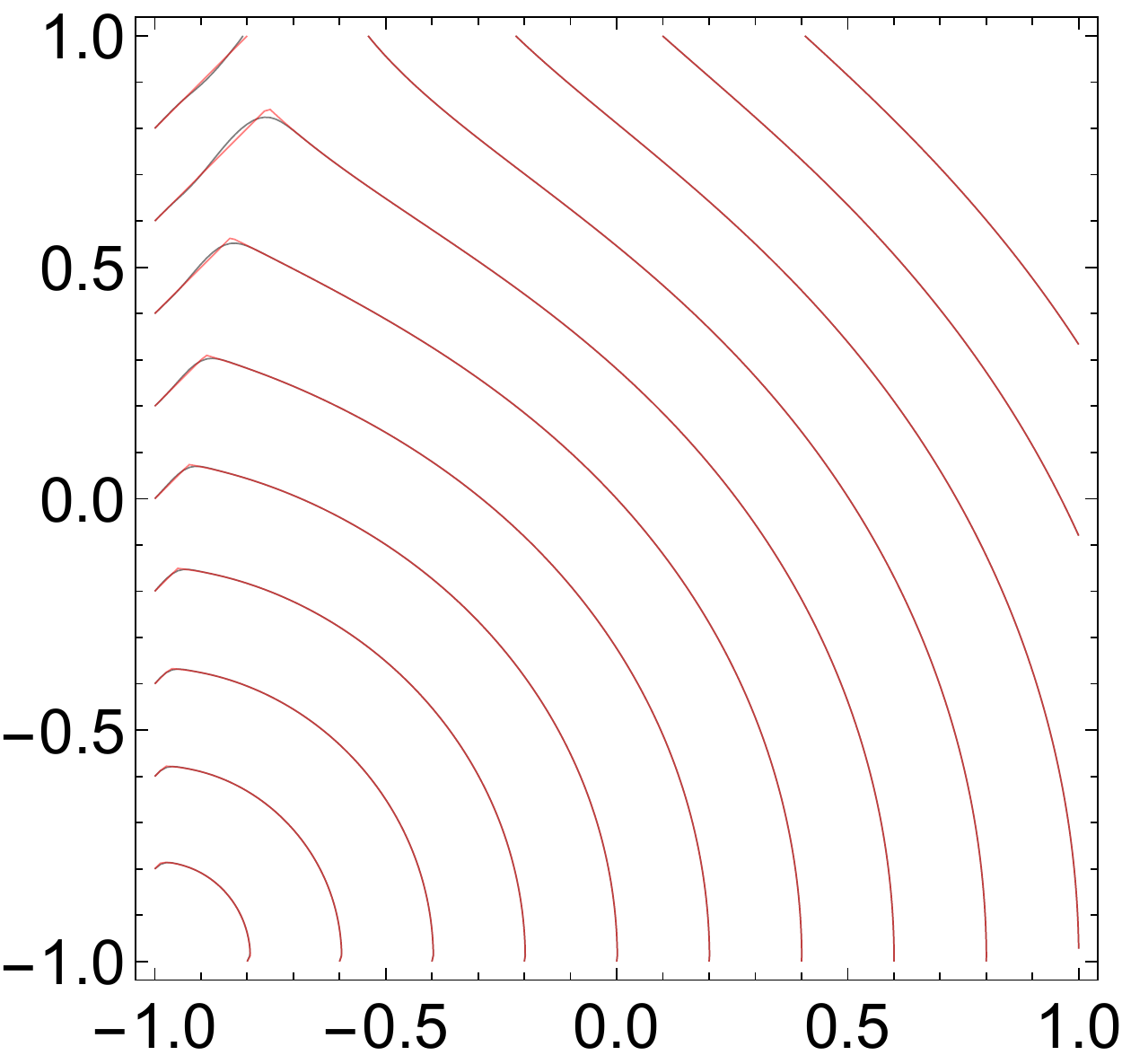}
\includegraphics[width=0.24\columnwidth]{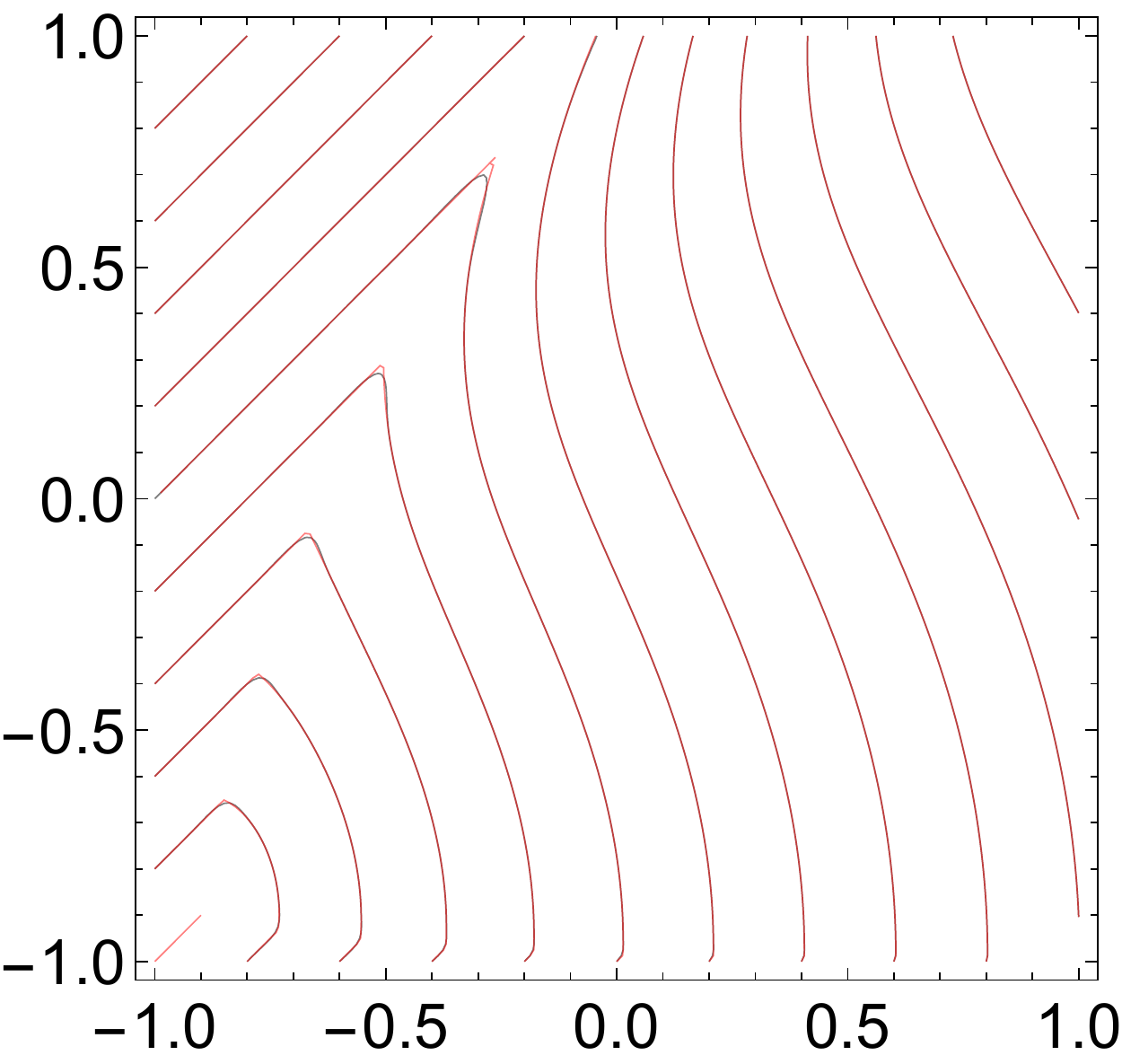}
\includegraphics[width=0.24\columnwidth]{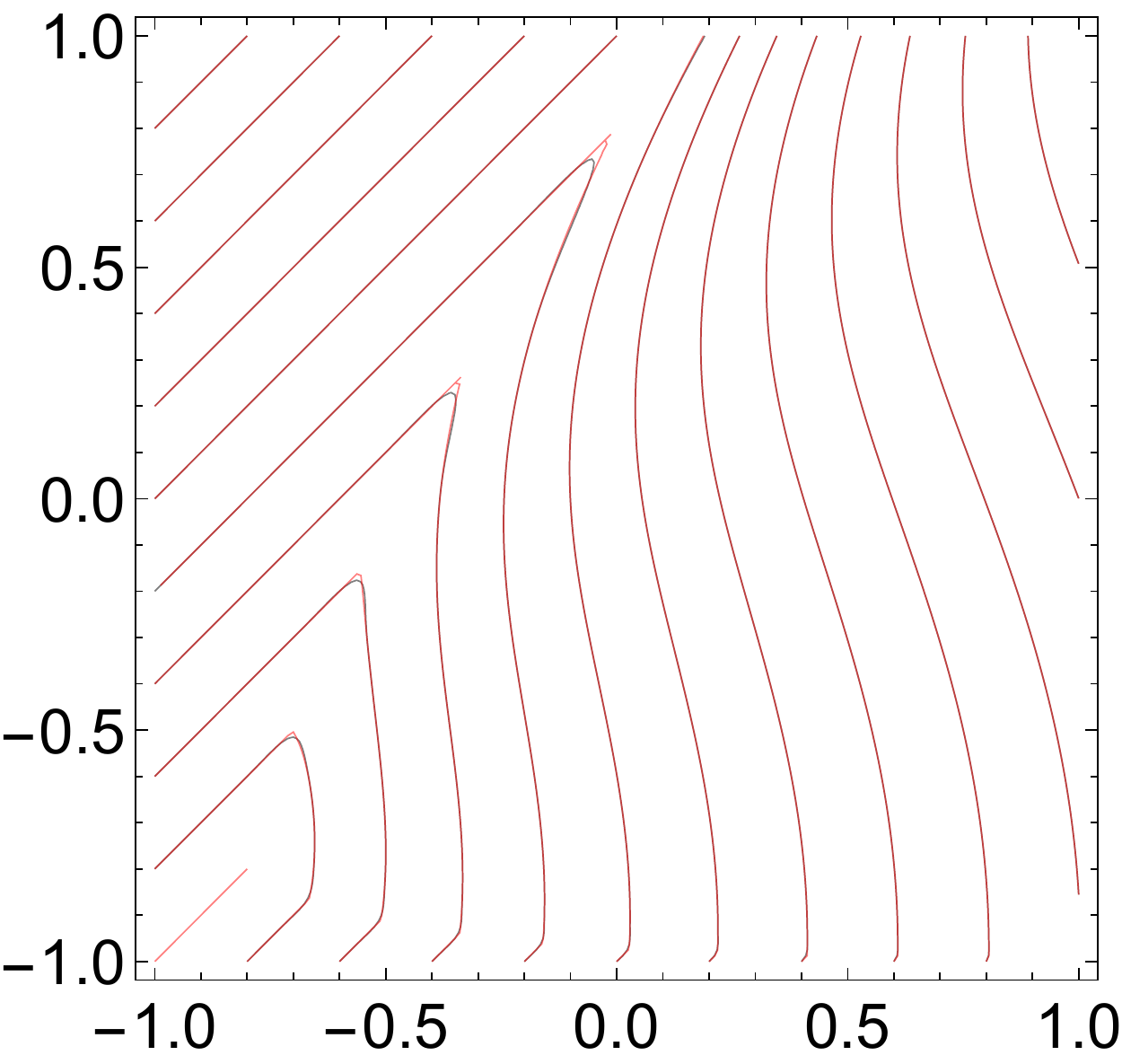}
\includegraphics[width=0.24\columnwidth]{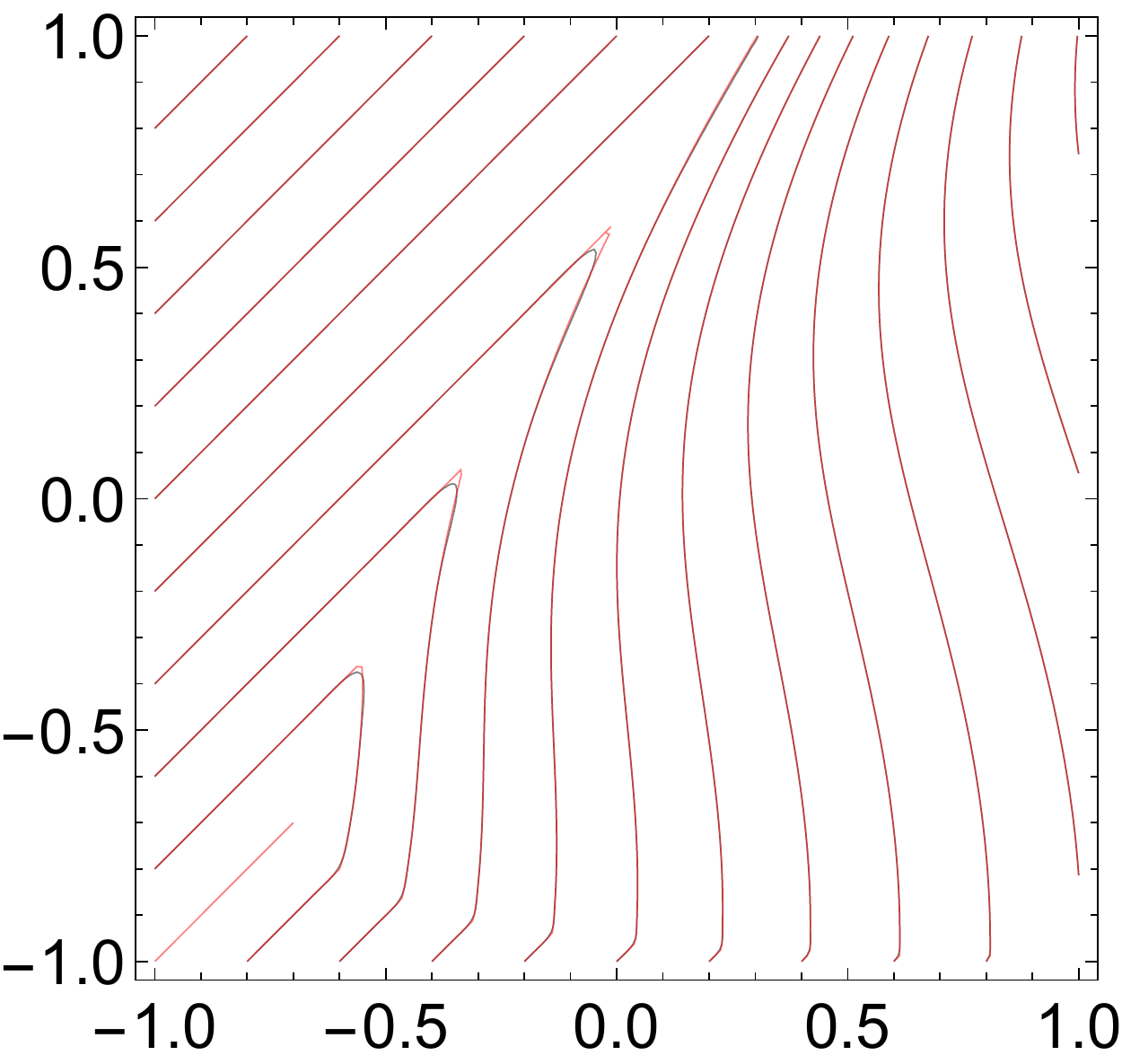}\\[1ex]
\includegraphics[width=0.24\columnwidth]{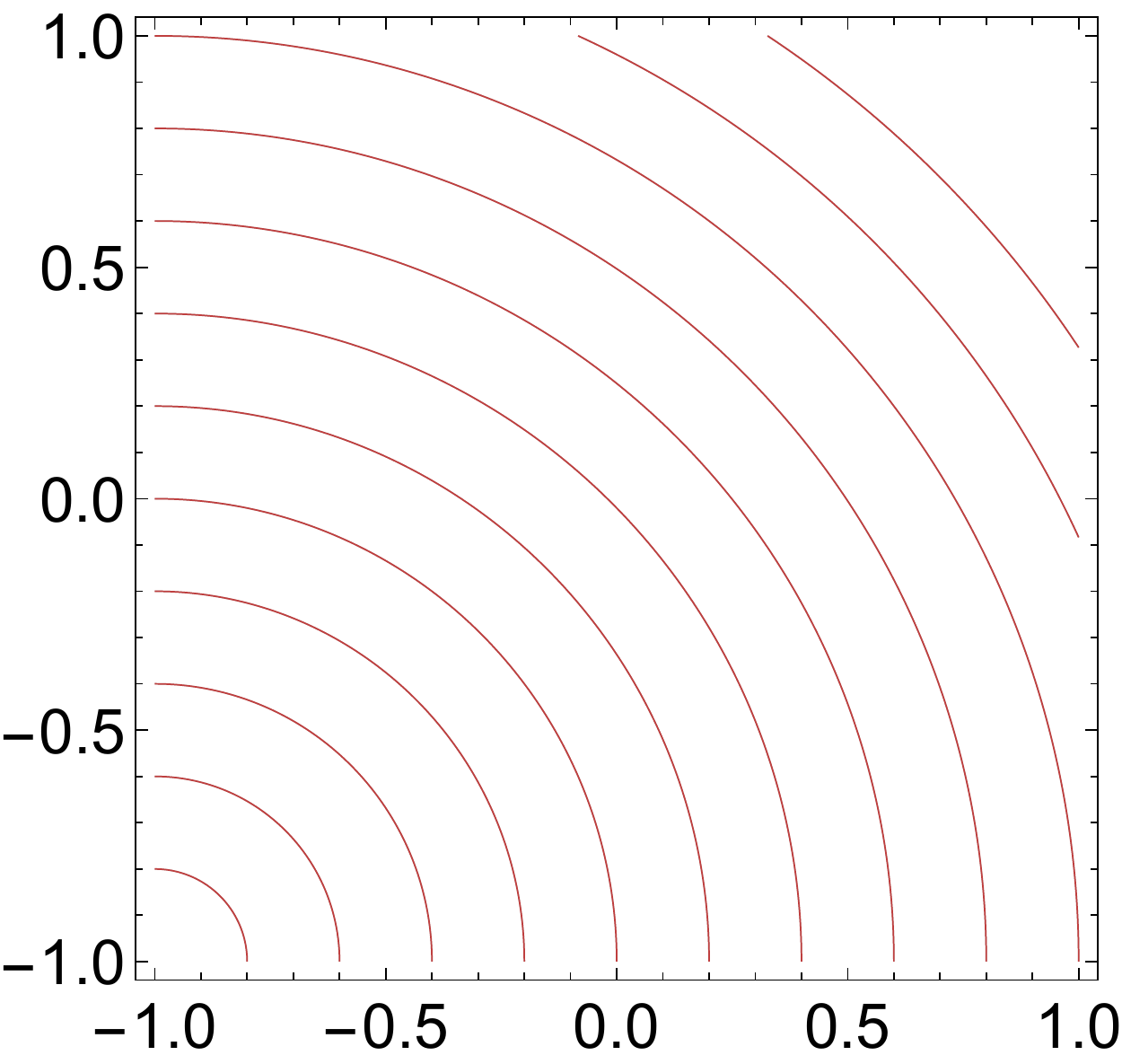}
\includegraphics[width=0.24\columnwidth]{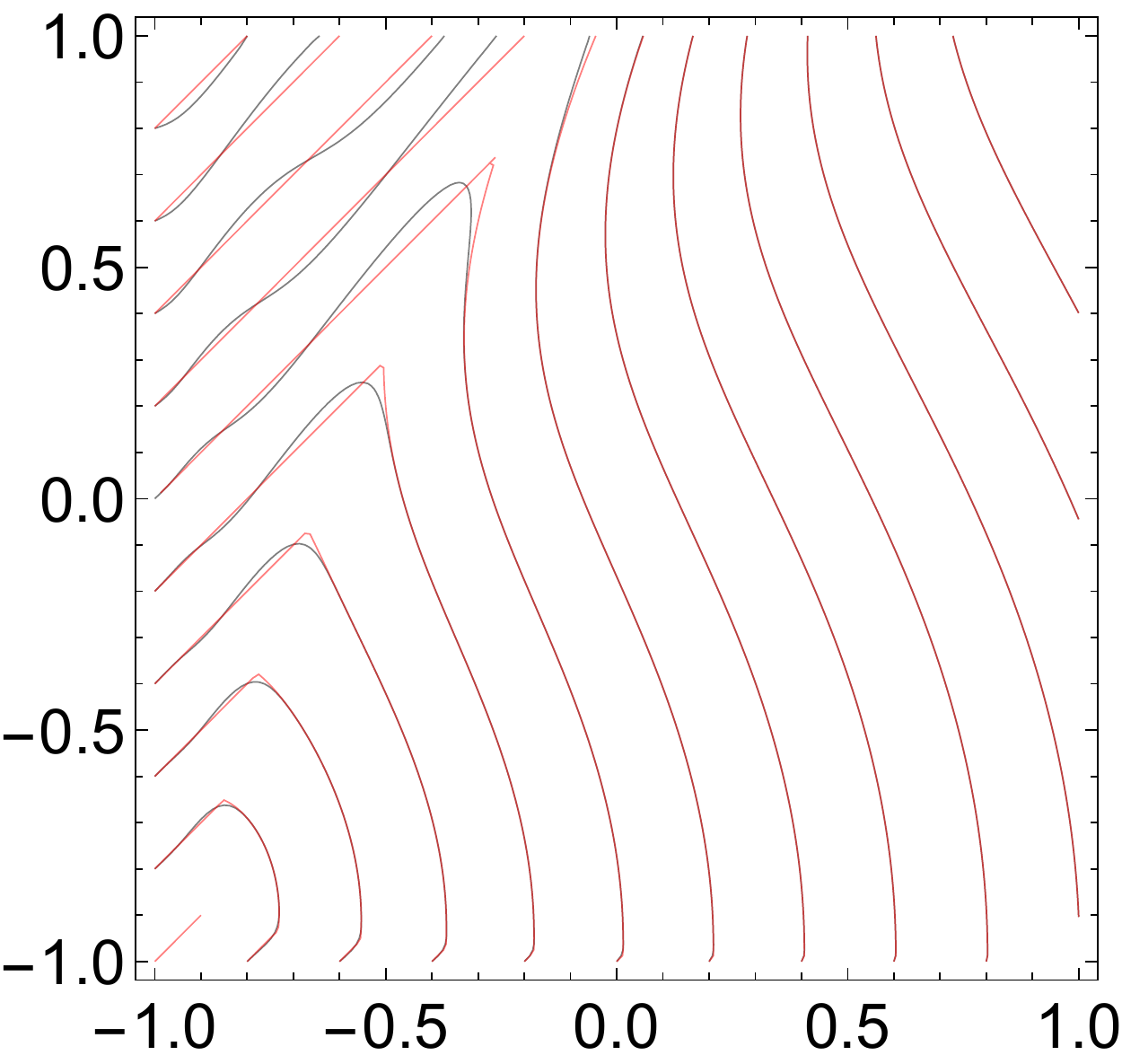}
\includegraphics[width=0.24\columnwidth]{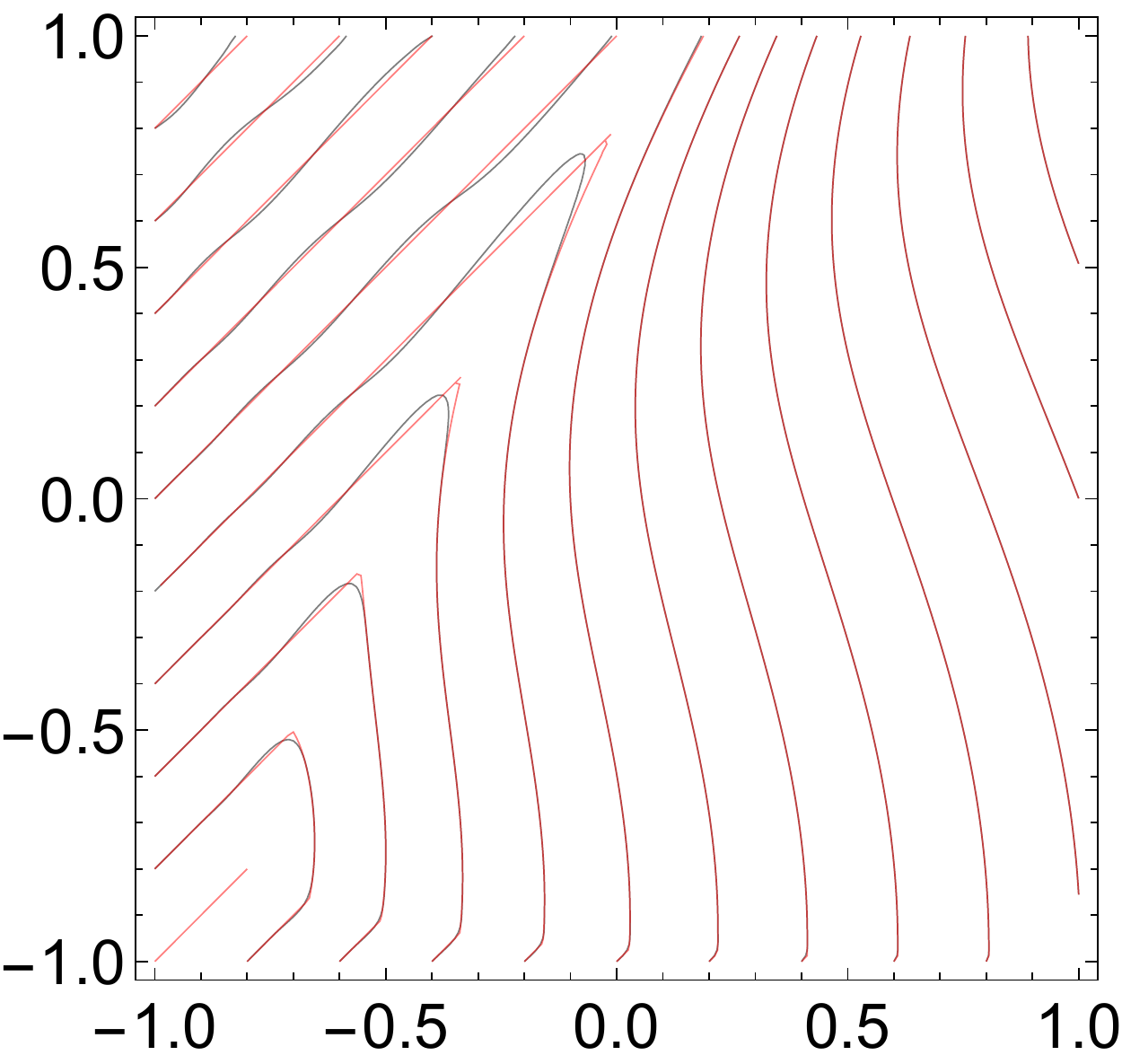}
\includegraphics[width=0.24\columnwidth]{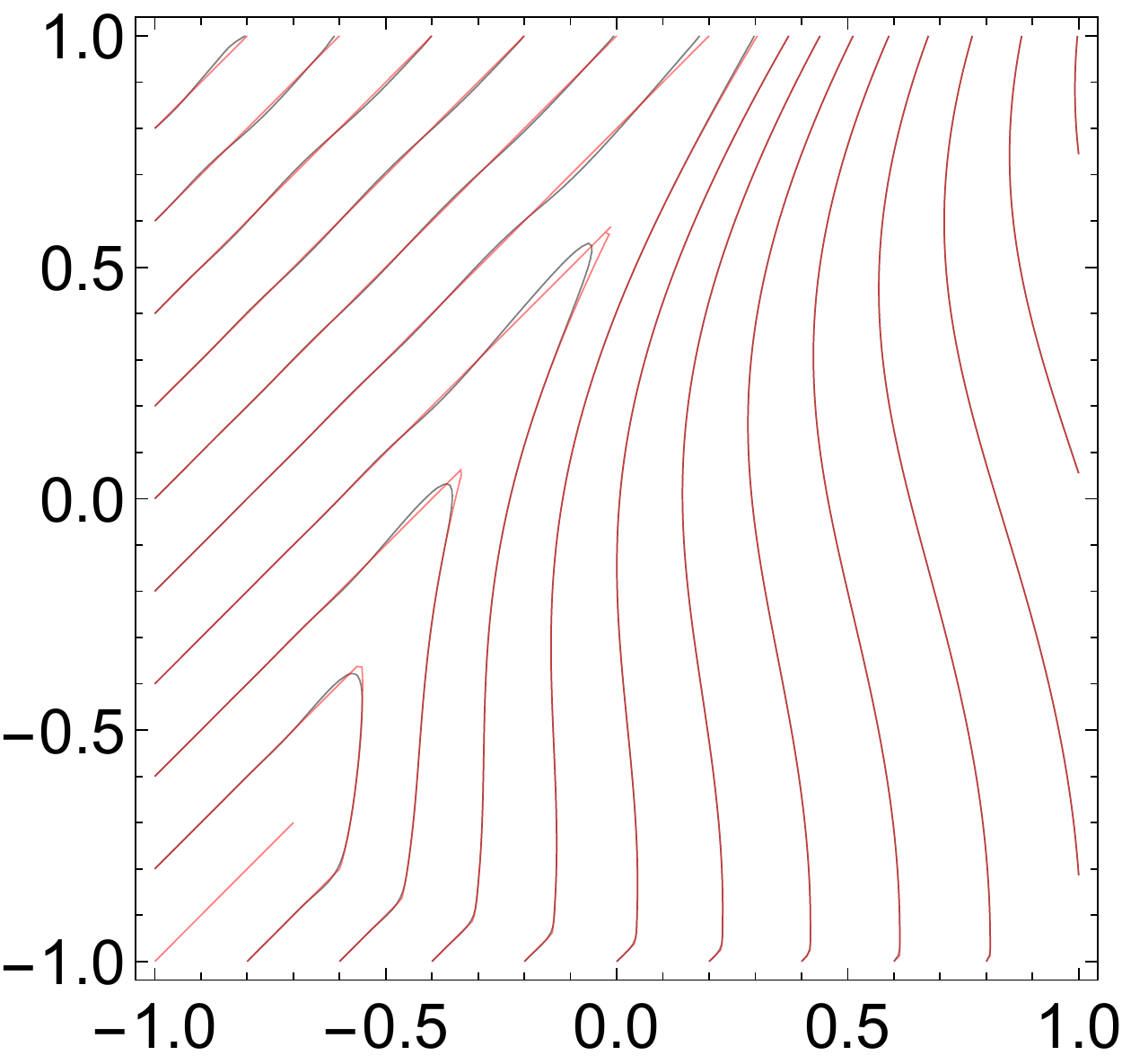}
\caption{Numerical solutions for the example with the exponential velocity for the finest grid $M=160$ and the times $t=0.1, 0.2, 0.3$ (from the second to the fourth column) using  (\ref{si2d}) with (\ref{k0}). The first row is obtained with $N=160$ when the maximal Courant number is $10.9$ (the first picture is the solution at $t=0.01$), the second row with $N=16$ with the maximal Courant number being $109$ (the first picture is the initial function $u^0$). 
The red contour lines represent the exact solution, the black contour lines represent the numerical solutions for values $0.0$, $0.2$ up to $2.4$. Note that no instabilities occur in the numerical solutions.
}
\label{figexp}
\end{center}
\end{figure}

We compute the example with the scheme (\ref{si2d}) using (\ref{kp}) - (\ref{k0}) and with the CTU scheme (\ref{ctusi2d1}) using (\ref{k3}) for a medium Courant number $10.9$ and for a large Courant number $109$. The results are summarized in Table \ref{tabexp} and in Figure \ref{figexp}. Note that no instabilities occur in the numerical solutions and the accuracy for the large Courant number is still acceptable. Only one sweep was used to solve the linear algebraic systems.

\begin{table}[h!]
\begin{center}
\begin{tabular}{||c||c|c||c|c||c|c||c|c||}
\hline
$M$ & $E_{10.9}$ & $E_{109.}$  & $E_{10.9}$ & $E_{109.}$ & $E_{10.9}$ & $E_{109.}$ & $E_{10.9}$ & $E_{109.}$ \\
\hline
40  & 33.5 & 99.7 & 24.7 & 103.  & 12.2 & 97.2 & 11.8 & 106. \\
80  & 13.8 & 44.7 & 10.1 & 45.4 & 4.29 & 44.1 & 3.92 & 51.0  \\
160& 5.67 & 19.8 & 4.05 & 18.7 & 1.55 & 18.9 & 1.34 & 24.6  \\
\hline
\end{tabular}
\caption{The error (\ref{error}) (multiplied by $10^3$) for the example with the exponential velocity for the maximal Courant numbers $10.9$ and $109$ using (\ref{si2d}) with (\ref{kp})  (the $2^{nd}$ - $3^{rd}$ columns), (\ref{km}) (the $4^{th}$ - $5^{th}$ ones), (\ref{k0}) (the $6^{th}$ - $7^{th}$ ones), and (\ref{ctusi2d1}) with  (\ref{k3}) . Note that the exact solution is not smooth, therefore the EOCs of the schemes are approaching the $1^{st}$ order accuracy from above. }
\label{tabexp}
\end{center}
\end{table}

\subsection{Two benchmark examples}

In the last section on numerical experiments we present two standard benchmark examples for the tracking of moving interfaces - the rotation of Zalesak's disc \cite{for05,fm07,Herrmann2008a,Kees2011,Kim2011,Wang2012,Samuel2014,Starinshak2014}
and the single vortex example \cite{for05,Min2007,fm07,Herrmann2008a,Kees2011,Kim2011,Wang2012,Samuel2014,Starinshak2014}
Note that in a standard setting of these examples there is no need to use a implicit or a semi-implicit scheme, nevertheless we compute these benchmark examples to show a satisfactory accuracy of the proposed semi-implicit $\kappa$-scheme.

In the first example the initial function $u^0(x,y)$ is a signed distance function to a slotted circle of diameter $0.6$ with the center in $(0,0.5)$ and with a cut obtained by a rectangle of the lengths $0.1 \times 0.5$ with the bottom corners at $(-0.05,y_r)$ and $(0.05,y_r)$ and $y_r=0.5-\sqrt{0.3^2 - 0.05^2)}$. The velocity is given in (\ref{rotation}), so the initial profile of $u^0$ returns after one rotation to its origin position at $t=1$.

\begin{figure}[h!]
\begin{center}
\includegraphics[width=0.32\columnwidth]{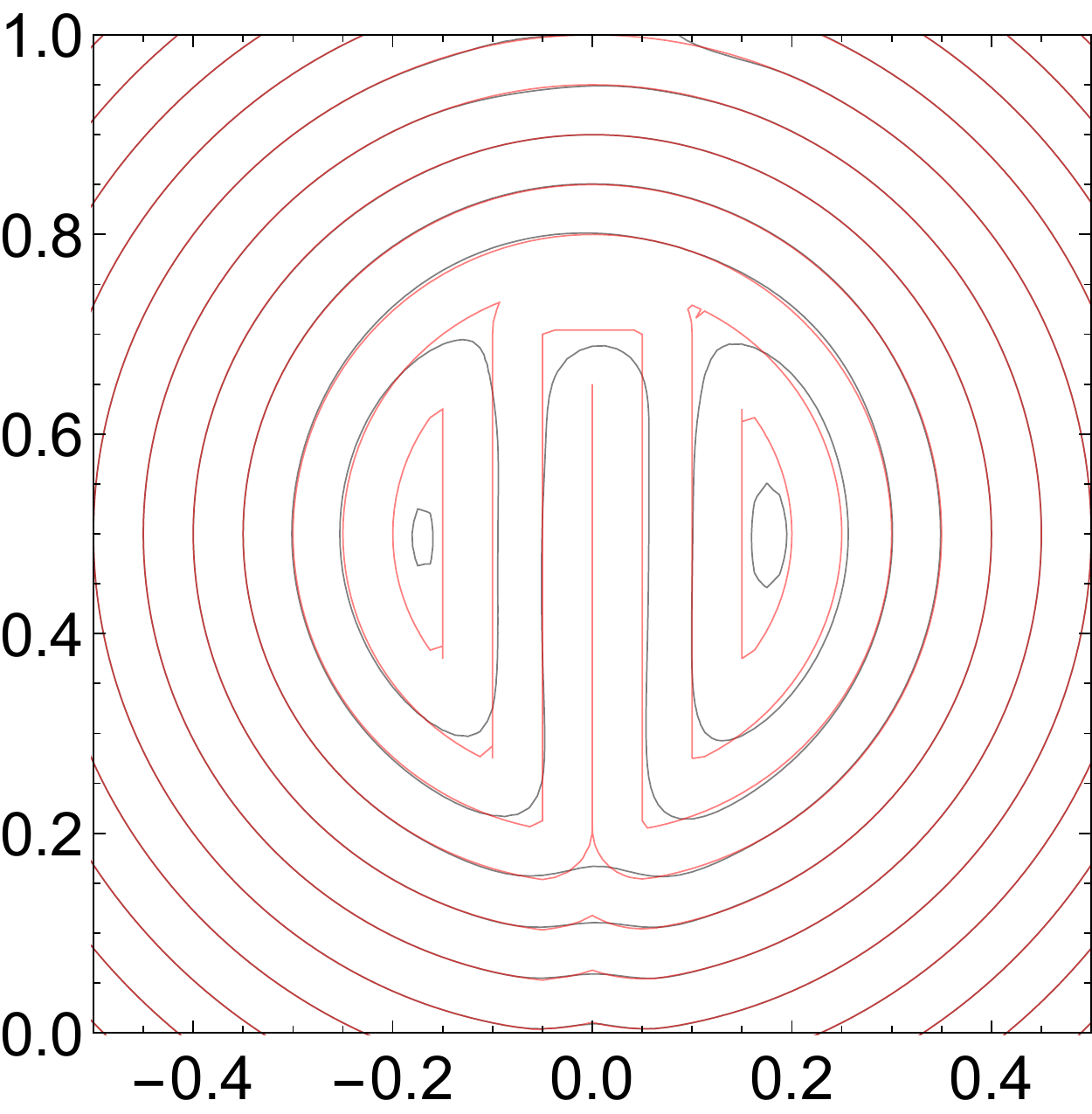}
\includegraphics[width=0.32\columnwidth]{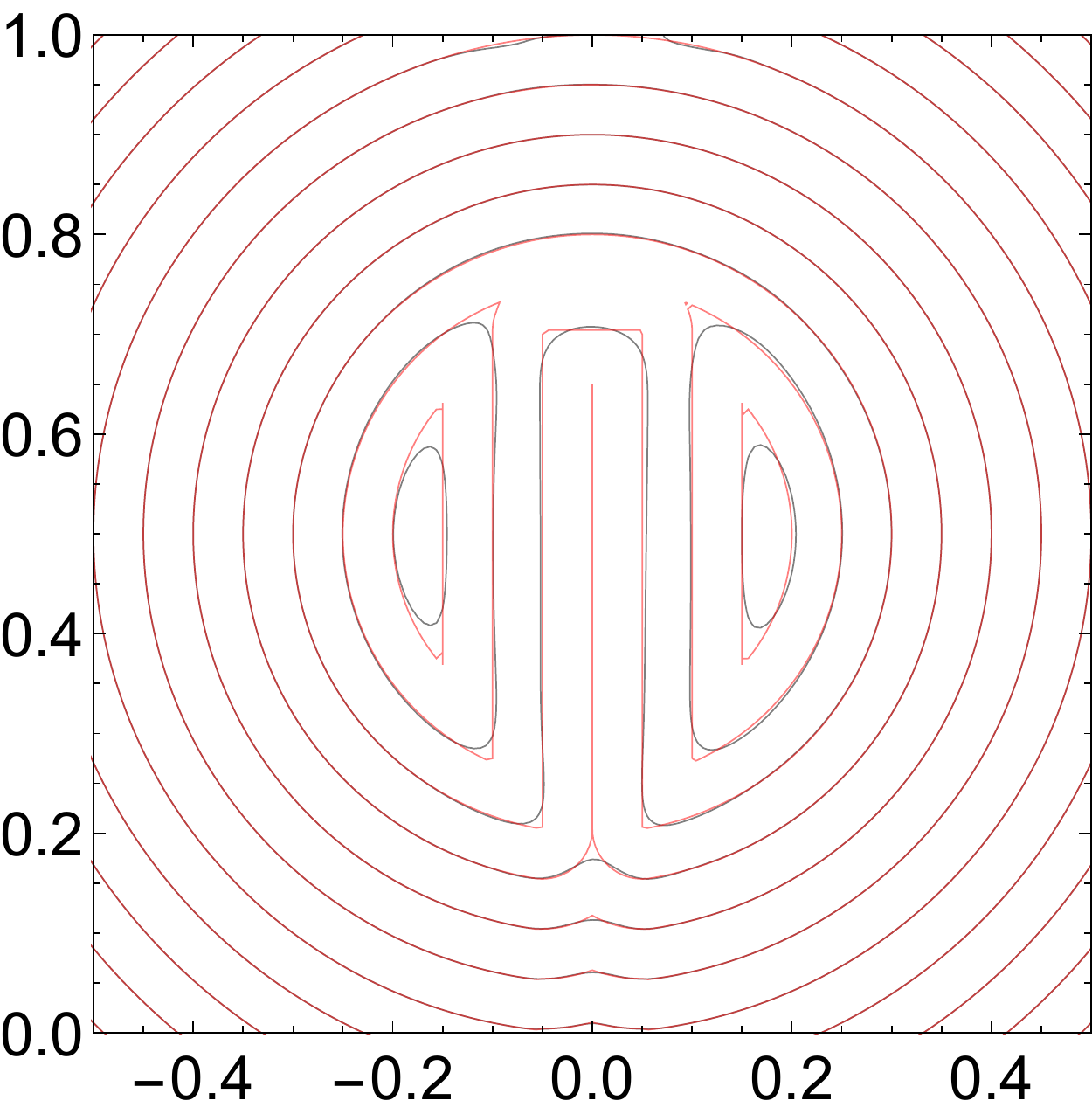}
\includegraphics[width=0.34\columnwidth]{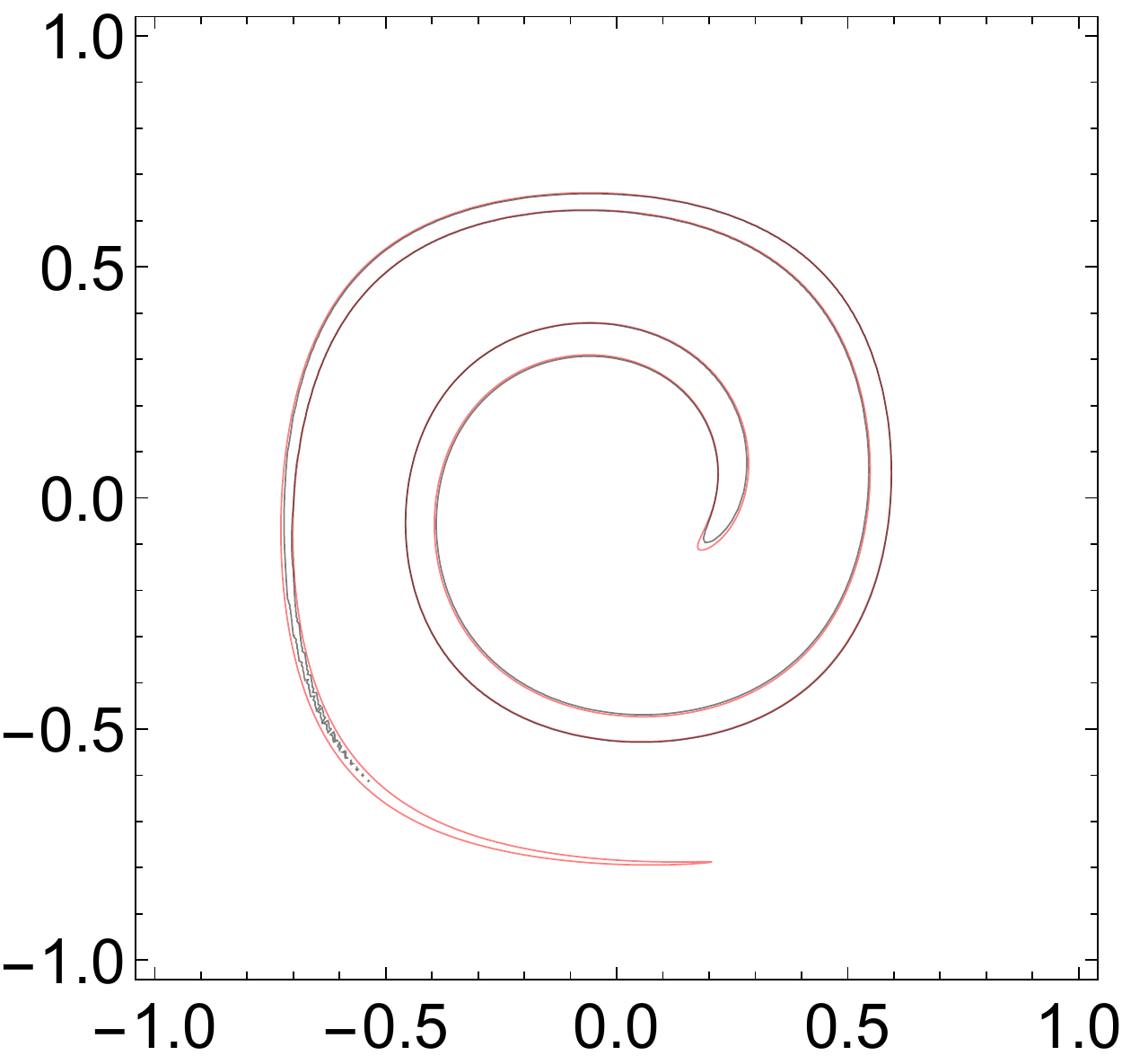}
\caption{Numerical solutions for the Zalesak's disc after one rotation computed using the scheme (\ref{si2d}) with (\ref{k3}) for the grids with $M=160$ (the first picture) and $M=320$ (the second picture). The red contour lines represent the exact solution, the black contour lines represent the numerical solutions for values $-0.05$, $0.0$ up to $0.4$. The third picture compares the numerical solutions for the single vortex example at $t=2.5$, the red contour line is the zero level set obtained with $M=1280$, the black contour line is the zero level set obtained with $M=160$ using always the scheme (\ref{si2d}) with (\ref{k3}).
} 
\label{figz}
\end{center}
\end{figure}

The EOCs for the scheme (\ref{si2d}) with all four variants (\ref{kp}) - (\ref{k3}) are presented in Table \ref{tab03}, and a visual grid convergence study for the scheme (\ref{si2d} with (\ref{k3}) is given in Figure \ref{figz}.

Next we choose the single vortex example that is characterized by a  deformational flow. The velocity $\vec{V}=(V,W)$ is given by
\begin{eqnarray*}
V(x,y)=-4 \sin^2(\pi(x+1)/2) \sin(\pi(y+1)/2) \cos(\pi(y+1)/2) \,,\\
W(x,y)=4 \sin^2(\pi(y+1)/2) \sin(\pi(x+1)/2) \cos(\pi(x+1)/2) \,.
\end{eqnarray*}
for $t \in [0,2.5]$. The initial function $u^0(x,y)$ is a signed distance function to a circle with the center at $(0,0.5)$ and the radius $0.3$. As $\vec{V} \equiv 0$ at the boundary $\partial D$, we fix the values of $u(x,y,t)=u^0(x,y)$ for $(x,y) \in \partial D$ and $t \in [0,2.5]$. 

We present the results at $t=2.5$ where the largest deformation of initial function can be observed. As the exact solution at this time is not available, we compute a reference numerical solution $\tilde U_{i j}^N$ obtained with  $N=5 {\cal M}/4$ and ${\cal M}=1280$ and compare it with numerical solutions for $M=80,160,320$ at $N=5 M/4$ using 
\begin{equation}
\label{l1e}
e = \frac{4}{M^2} \sum_{i,j=0}^M  |U_{i j }^{5M/4} - \tilde U_{i {\cal M}/M \,\, j {\cal M}/M }^{5{\cal M}/4}| \,.
\end{equation}
The chosen time step corresponds to the maximal Courant number equal to 2.
The results are summarized in Table \ref{tab03} and a visual comparison of the zero level set for the numerical solutions obtained with $M=1280$ and $M=160$ is given in Figure \ref{figz}. 

\begin{table}[h!]
\begin{center}
\begin{tabular}{||c||c|c||c|c||c|c||c|c||}
\hline
$M$ & $E_{z}$ & $E_{s}$ & $E_{z}$ & $E_{s}$ & $E_{z}$ & $E_{s}$ & $E_{z}$ & $E_{s}$  \\
\hline
40  & 57.3 &        & 32.9 &        & 15.3 &          & 15.2 &      \\
80  & 33.2 & 186. & 22.1 & 160. & 7.45 & 81.4 & 7.09 & 80.0\\
160& 15.6 & 74.3 & 12.3 & 65.0 & 3.54 & 30.1 & 3.21 & 28.2 \\
320&        & 25.6 &        & 25.4 &        & 8.50 &        & 7.70 \\
\hline
\end{tabular}
\caption{The error (\ref{l1e}) (multiplied by $10^3$)  for the numerical solutions of Zalesak's disc (the columns with $E_{z}$) and the single vortex (the columns with $E_{s}$) that are solved using (\ref{si2d}) with (\ref{kp}) (the $2^{nd}$ and $3^{rd}$ column), (\ref{km}), (\ref{k0}), and (\ref{k3}) (the $8^{th}$ and $9^{th}$ columns). Note that $N=5 M/4$ that corresponds to the maximal Courant number for the rotation of Zalesak's disc equal to $2.5$ and for the single vortex example equal to $2$.}
\label{tab03}
\end{center}
\end{table}

\begin{figure}[h!]
\begin{center}
\includegraphics[width=0.3\columnwidth]{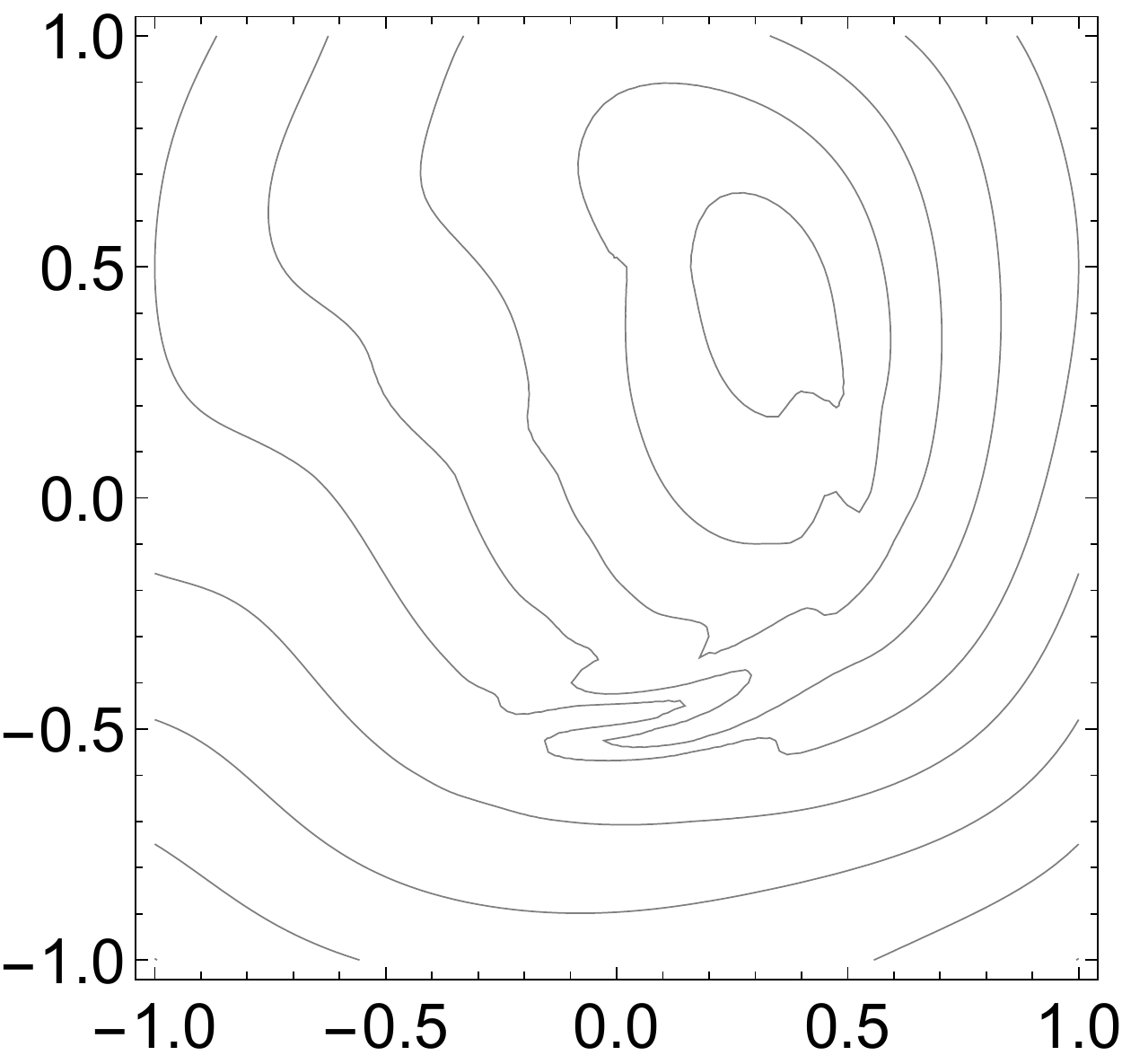}
\includegraphics[width=0.3\columnwidth]{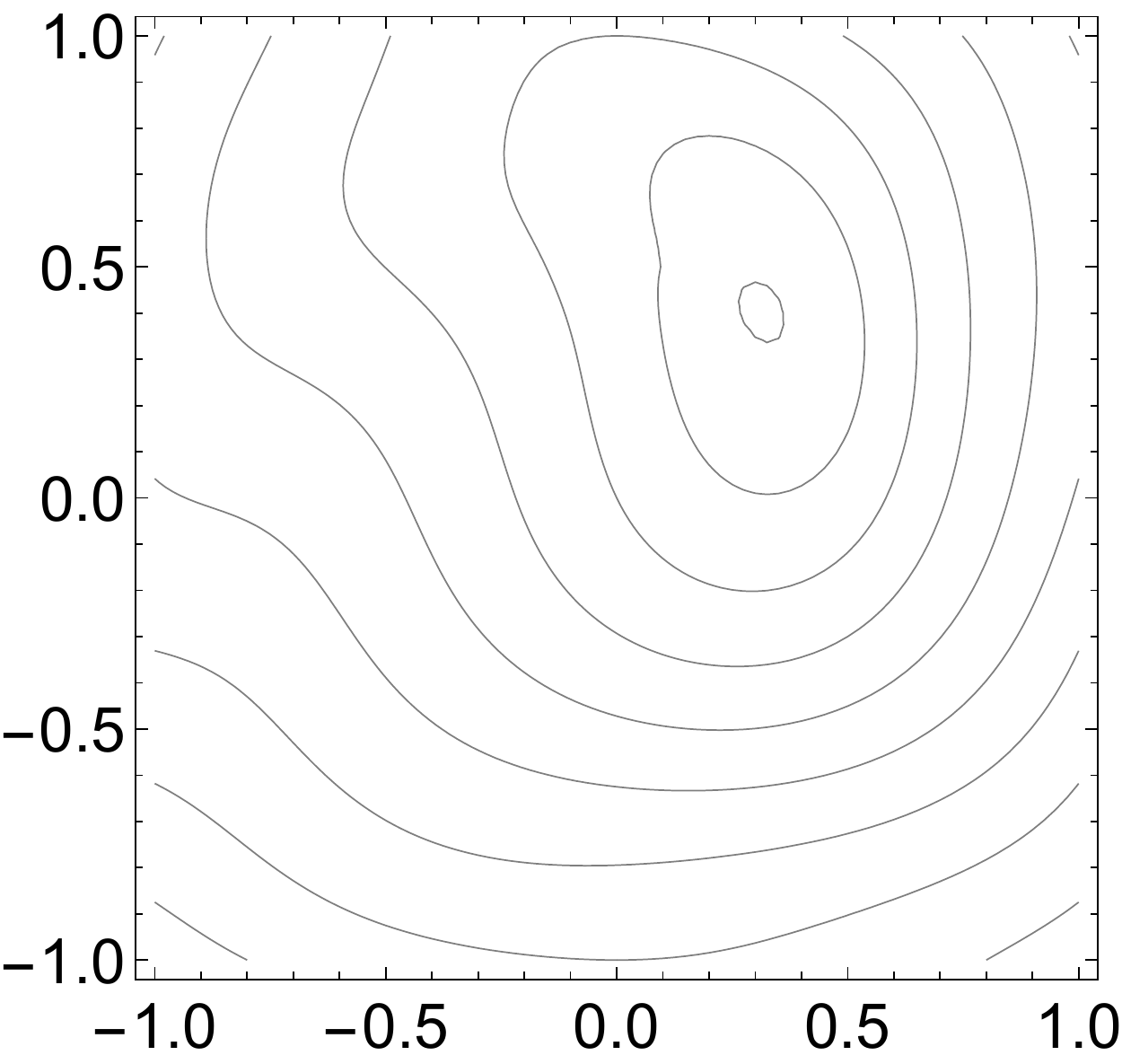}
\includegraphics[width=0.3\columnwidth]{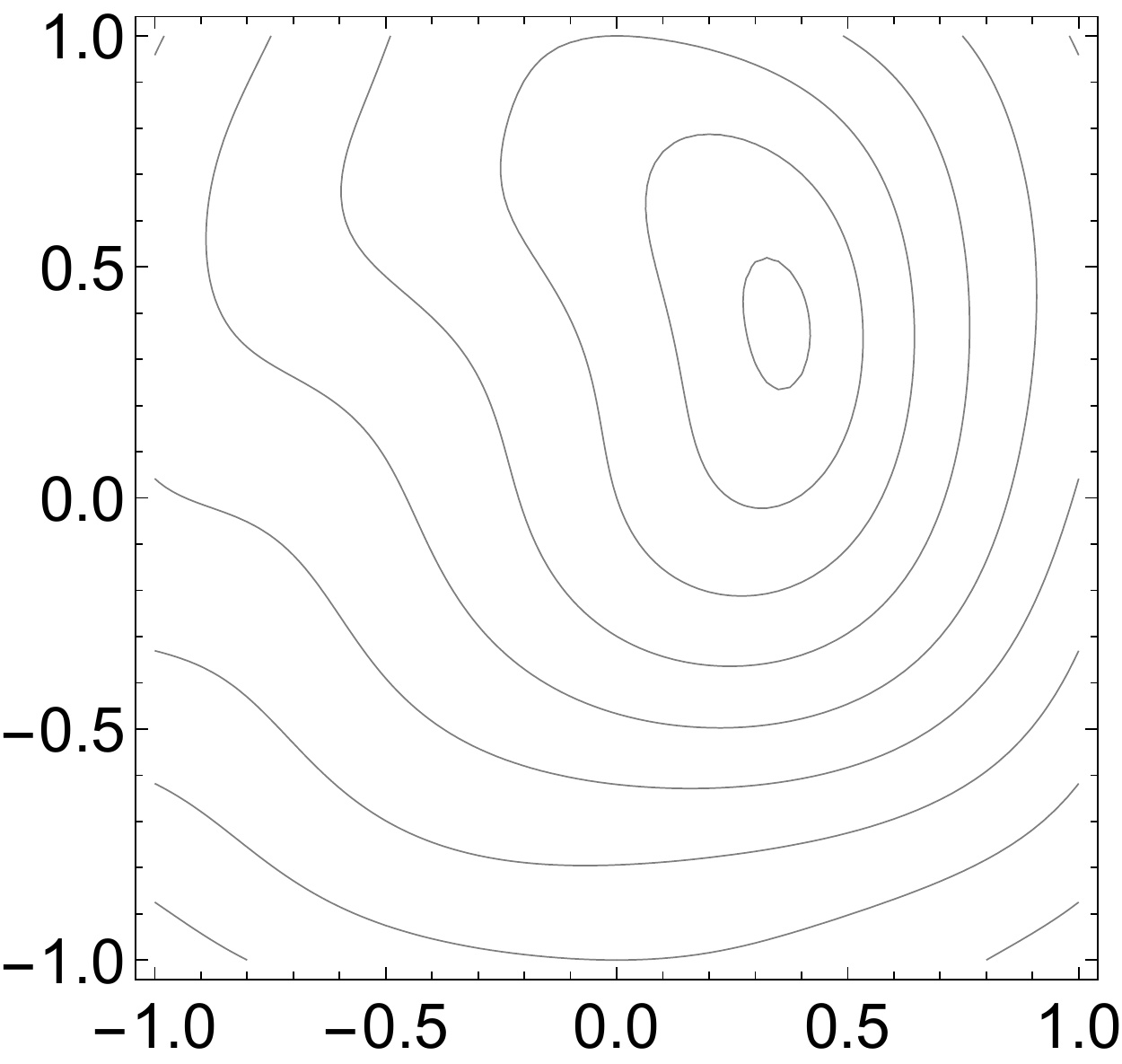}
\caption{Numerical solutions at $t=0.2$ with contour lines for values $0.1$, $0.3$ up to $1.7$ for the single vortex example for the coarse grid $M=80$. The first picture is obtained using one time step (the maximal Courant number equals $16$) with the scheme (\ref{si2d}) with (\ref{k3}) and no CTU extension when instabilities can be observed. The second picture is obtained using the CTU extension (\ref{ctusi2d1}) with (\ref{k3}) taking one identical time step when no instabilities occur. For a comparison the last picture is obtained using $16$ uniform time steps that corresponds to the maximal Courant number being $1$ using (\ref{si2d}) with (\ref{k3}).
}
\label{fig04}
\end{center}
\end{figure}

Finally, we can illustrate for the single vortex example the stability results obtained with the numerical von Neumann stability analysis.
Firstly we compute the example with the scheme (\ref{si2d}) with no CTU extension using (\ref{k3}) applying one large time step such that the maximal Courant number equals $16$. In Figure \ref{fig04} we present the numerical result obtained after one sweep where some instabilities can be clearly observed. Computing the example with the identical one large time step but using the CTU extension (\ref{ctusi2d1}) with (\ref{k3}), no instabilities occur, see Figure \ref{fig04}. 

\section{Conclusions}
\label{sec-con}

In this paper the class of semi-implicit schemes on Cartesian grids for the numerical solution of linear advection equation is derived. The derivation follows the partial Lax-Wendroff (or Cauchy-Kowalewski) procedure to replace the time derivatives of the exact solution in Taylor series. Opposite to the full form of this procedure when only the space derivatives are used in the replacement, the partial procedure exploits also the mixed time-space derivatives. 

The one-dimensional semi-implicit $\kappa$-scheme (\ref{si1d})  is $2^{nd}$ order accurate with the unconditional stability for the variable velocity case and for all considered values of $\kappa$. The analogous fully implicit $\kappa$-scheme (\ref{impl1d}) is unconditionally stable only for a limited range of $\kappa$ and only for the locally frozen values of advection velocity. 

The dimension by dimension extension (\ref{si2d}) of one-dimensional semi-implicit $\kappa$-scheme to Cartesian grids in several dimensions gives the $2^{nd}$ order accurate scheme. The analogous extension of fully implicit $\kappa$-scheme results only in a first order accurate scheme. 

We derive the Corner Transport Upwind extension of two-dimensional semi-implicit $\kappa$-scheme. The semi-implicit $\kappa$-scheme (\ref{ctusi2d1}) with the variable choice (\ref{kappasi2d}) of $\kappa$ parameters has all desired properties that are considered in this paper - the scheme is $2^{nd}$ order accurate for the variable advection velocity and unconditionally stable according to the numerical von Neumann stability analysis. Moreover it is $3^{rd}$ order accurate if the velocity is constant.  The chosen examples in numerical experiments confirm these properties.

\end{document}